\documentclass{article}
\RequirePackage{lineno}
\usepackage[textwidth=16cm]{geometry}

\usepackage[footnotesize,tt]{caption}
\usepackage{subcaption}
\usepackage{amsfonts}

\usepackage{pdfpages}
\usepackage{authblk}
\usepackage{graphics}
\usepackage{mathrsfs}

\usepackage{graphicx}

\patchcmd{\thebibliography}{\chapter*}{\section*}{}{}
\usepackage{multirow}

\usepackage[most]{tcolorbox}
\usepackage{caption}

\usepackage{siunitx}
\usepackage{amsthm}
\numberwithin{equation}{section}
\usepackage{amsmath,bm}
\usepackage{mathtools}
\usepackage[hidelinks]{hyperref}
\usepackage{multicol}
\usepackage{tikz}
\usepackage{float}
\usepackage{amsmath,amssymb,amsthm}
\usepackage{bm}
\usepackage{mathtools}

\usetikzlibrary{fadings}
\usepackage[utf8]{inputenc}
\usetikzlibrary{fit,positioning}
\usetikzlibrary{positioning,calc}

\usepackage{amssymb}

\usepackage{amsthm}
\usepackage[hidelinks]{hyperref}
\usepackage{listings}
\usetikzlibrary{positioning,decorations.pathreplacing,shapes}
\usepackage[english]{babel}
\usepackage{microtype}
\usepackage{xcolor}
\usepackage{soul}

\tikzset{
  pivot/.style={
    draw,
    regular polygon,
    regular polygon sides = 3,
    fill = cyan,
    node distance = 1cm,
    minimum height = 2cm,
    at = {(0,0)}
  }
}

\usepackage{colortbl}
\usepackage{amsthm}

\newcommand{\Om}{\Omega}
\newcommand{\dOm}{\partial\Omega}
\newcommand{\dd}{\,\mathrm{d}}
\newcommand{\bn}{\bm{n}}
\newcommand{\bu}{\bm{u}}
\newcommand{\bx}{\bm{x}}
\newcommand{\bsig}{\bm{\sigma}}
\newcommand{\bveps}{\bm{\varepsilon}}
\newcommand{\divv}{\operatorname{div}}
\newcommand{\tr}{\operatorname{tr}}
\newcommand{\grad}{\nabla}

\newtheorem{theorem}{Theorem}[section]

\newtheorem{remark}{Remark}

\title{A 2D Free Boundary Model for Atherosclerosis with Diffusion and Elasticity}

\author[1]{G. Abi Younes}
\author[1]{N. El Khatib\thanks{Corresponding author: nader.elkhatib@lau.edu.lb}}
\author[2]{N. Bessonov}
\author[3,4]{V. Volpert}

\affil[1]{\textit{\small Lebanese American University, Department of Computer Science and Mathematics, Byblos, Lebanon}}
\affil[2]{\textit{\small Institute of Mechanical Engineering Problems, 199178, Saint Petersburg, Russia}}
\affil[3]{\textit{\small Institut Camille Jordan, UMR 5208 CNRS, University Lyon 1, 69622 Villeurbanne, France}}
\affil[4]{\textit{\small Institute of Fundamental Technological Research, Polish Academy of Sciences 5B Pawinskiego St., 02-106 Warsaw, Poland}}

\date{}

\begin{document}

\maketitle

\begin{abstract}
We propose a coupled elasticity--diffusion model on a deformable domain to describe the interaction between diffusion and mechanical deformation. We first study the two-dimensional formulation and investigate its numerical behavior. A quasi-stationary formulation is then derived, and a linear relation between the total mass of the diffusing species and the volume change of the domain is established. Next, a one-dimensional reduction is introduced as a particular case of the two-dimensional model to provide additional analytical insight into the coupled dynamics. Finally, the proposed framework is applied to a model of atherosclerosis, where numerical simulations illustrate the influence of the model parameters on inflammatory wave propagation and arterial wall deformation.

\end{abstract}

\section{Introduction}
Partial differential equations play a central role in modeling physical and biological systems. In particular, reaction–diffusion equations describe the transport and interaction of chemical species, while elasticity equations govern the deformation of the media. In many applications, these processes are coupled: chemical activity affects the mechanical deformation of the medium, while the deformation influences the distribution of chemical substances. When the domain itself evolves over time, the governing equations must be solved together with the domain evolution, leading to additional analytical and numerical challenges. These interactions are especially important in biological tissues and are encountered in diseases such as atherosclerosis. This motivates the study of mathematical models for atherosclerosis that combine reaction, diffusion, elasticity, and domain evolution.

Atherosclerosis is a chronic inflammatory disease characterized by the accumulation of lipid-rich plaques within the arterial wall \cite{libby2002}. It is initiated by endothelial dysfunction and lipid infiltration, which trigger an inflammatory response involving immune cells and cytokines \cite{libby2009}. As the disease progresses, plaques may narrow or rupture arteries, leading to serious cardiovascular events such as ischemic heart disease and stroke \cite{libby2009}. Cardiovascular diseases are the leading cause of death worldwide, with atherosclerosis being a major underlying cause \cite{who2025}.

Low-density lipoprotein (LDL) cholesterol is a major risk factor in the development of atherosclerosis \cite{libby2002}. LDL particles penetrate the arterial wall and undergo oxidative modification, forming oxidized LDL (ox-LDL) \cite{steinberg2009}. Oxidized LDL triggers an inflammatory response by activating endothelial cells and promoting the release of chemoattractants \cite{hansson2005}. This process recruits circulating monocytes, which migrate into the arterial wall and differentiate into macrophages, further amplifying inflammation \cite{libby2002}.

Macrophages internalize ox-LDL through scavenger receptors and transform into lipid-laden foam cells, a hallmark of early atherosclerotic lesions \cite{glass2001}. Foam cells contribute to the formation of fatty streaks, and their persistence sustains chronic inflammation. They secrete pro-inflammatory cytokines, including tumor necrosis factor (TNF-$\alpha$) and interleukin-1 (IL-1), which activate endothelial cells, recruit additional monocytes, and amplify the inflammatory response \cite{libby2002,hansson2005}.

This inflammatory response is partially regulated by anti-inflammatory cytokines, such as interleukin-10 (IL-10), which inhibit the production of pro-inflammatory cytokines \cite{moore2001}. Meanwhile, smooth muscle cells migrate from the media to the intima, forming a fibrous cap over the lipid core \cite{ross1999}. Although the fibrous cap stabilizes the plaque, it also alters the arterial geometry and blood flow, which may promote plaque rupture or thrombus formation \cite{libby2002,virmani2000}.

The arterial wall exhibits elastic properties that are essential for vascular function. Under physiological conditions, it undergoes cyclic, reversible deformations in response to pulsatile blood flow \cite{humphrey2002}. This fluid--structure interaction contributes to blood pressure regulation and efficient circulation \cite{nichols2011}. In atherosclerosis, however, lipid accumulation and inflammation alter the composition of the arterial wall \cite{libby2002}, while smooth muscle cell proliferation and fibrous cap formation modify its mechanical structure \cite{stary1995}. These changes affect the mechanical properties of the wall and promote structural deformation \cite{humphrey2002b}, altering the vessel geometry and increasing the risk of plaque instability and rupture.

In this work, we consider a coupled system consisting of two-dimensional elasticity equations and a reaction--diffusion equation, where the domain deformation is governed by the mechanical response of the medium. We apply this model to atherosclerosis to study the interaction between inflammation and arterial wall deformation during disease progression.

\section{Literature Review}

The mathematical modeling of atherosclerosis has attracted considerable
attention over the past two decades, driven by the need to understand the
complex interplay between biochemical processes, immune-cell dynamics, and
mechanical deformation of the arterial wall.  The literature on this subject
spans several distinct but interconnected directions: reaction--diffusion
models of inflammatory dynamics, free-boundary and moving-domain formulations
of plaque growth, mechanical models of arterial-wall behavior, and coupled
chemo-mechanical systems.  The present work draws on all of these threads and
proposes a unified two-dimensional framework that integrates linear elasticity
with diffusion and reaction on an evolving domain.

\medskip
\noindent\textbf{Mathematical models of atherosclerosis.}
Mathematical models of atherosclerosis have addressed the biochemical processes underlying lesion initiation and progression. In this context, El
Khatib, G\'enieys, and Volpert~\cite{ElKhatib2007} introduced a
reaction--diffusion system capturing the inflammatory response in the arterial
intima, modelling the interactions between macrophages and cytokines and
demonstrating the existence of travelling-wave solutions that describe the
spatial propagation of inflammation.  This framework was subsequently extended
to incorporate LDL oxidation, monocyte recruitment, and feedback mechanisms
involving pro- and anti-inflammatory cytokines~\cite{elkhatib2012,
bulelzai2012}.  The role of oxidized LDL as the primary biochemical trigger was
analyzed rigorously in the context of threshold phenomena and bifurcation
theory, connecting the onset of a self-sustained inflammatory state to the
magnitude of the initial lipid infiltration.  Parallel modelling efforts by
Calvez, Ebde, Meunier, and Raoult~\cite{calvez2009} and by Chalmers, Cohen,
Terry, and Chambers~\cite{chalmers2015} investigated reduced ODE and PDE
systems to characterize steady states and the transition between healthy and
pathological regimes.  More recently, Parton et al.~\cite{parton2016}
incorporated stochastic components to account for inter-individual variability
in plaque progression.  The biological foundations underlying these models ---
including endothelial dysfunction, monocyte transmigration,
macrophage-to-foam-cell conversion, and cytokine signalling --- are thoroughly
documented in the reviews of Libby~\cite{libby2002}, Hansson~\cite{hansson2005},
Ross~\cite{ross1999}, and Libby, Ridker, and Hansson~\cite{libby2009}, which
collectively establish the immunological basis for the mathematical descriptions
adopted in the present paper.

\medskip
\noindent\textbf{Free-boundary problems in biological and vascular
modelling.}
The treatment of plaque growth as a free-boundary problem, in which the
evolving geometry of the arterial wall is itself an unknown coupled to the
underlying field equations, represents one of the most mathematically
challenging aspects of atherosclerosis modelling.  Free-boundary formulations
have a long history in mathematical biology, originating with models of tumour
growth and avascular spheroid dynamics.  The foundational work of Byrne and
Chaplain~\cite{byrne1995} and Greenspan~\cite{greenspan1976} established the
analytical framework for Stefan-type moving-boundary problems in biological
contexts, in which the normal velocity of the interface is determined by a local
balance condition.  In the context of atherosclerosis, Friedman and
Hao~\cite{friedman2015} proposed a free-boundary model for plaque growth in
which the domain boundary evolves according to the displacement induced by cell
accumulation, establishing well-posedness results for the associated PDE system.
Subsequent work by Friedman, Hao, and Hu~\cite{friedman2015b} analyzed the
stability of steady plaques and identified conditions under which the lesion
grows monotonically or reaches a finite equilibrium.  
In a related direction, models for the growth and remodelling of soft biological
tissues posed on evolving domains have been analyzed using the theory of
morpho-elasticity; see, e.g., Goriely and Ben Amar~\cite{goriely2005} and the
monograph of Goriely~\cite{goriely2017}, which provide the geometric and
mechanical underpinnings for domain evolution driven by growth tensors or
concentration-induced volumetric changes. Within this general morpho-elastic framework, the model of Goriely and
Ben Amar ~\cite{goriely2005} offers a representative example for studying the
interaction between growth and mechanics in biological tissues. Similar to the
present model, it investigates how growth-induced mechanical stresses can
generate tissue deformation and modify the geometry of the domain. However,
the two approaches differ significantly in their assumptions and objectives. Goriely and Ben Amar consider a three-dimensional spherical shell and assume radial symmetry, which reduces the mechanical problem to a one-dimensional radial configuration. Growth is prescribed through a morphoelastic growth tensor within a finite elasticity framework, and the main objective is to analyze the residual stresses generated by differential growth and their role in triggering buckling instabilities. In contrast, the present model is formulated in a two-dimensional moving domain without symmetry assumptions and couples quasi-static linear elasticity with a reaction--diffusion equation. The deformation is driven by the evolution of the concentration field rather than by a prescribed growth law, and the focus is placed on the coupled dynamics between diffusion, elasticity, and boundary evolution.\\

Classical Stefan-type formulations determine boundary motion through an interfacial flux condition. In the present setting, however, domain evolution is instead associated with the mechanical displacement induced by cellular accumulation within the tissue.\\

\medskip

Several subsequent works have explored the interaction between biochemical signaling, mechanics, and tissue growth on evolving domains. For example, Bressan and Lewicka ~\cite{BressanLewicka2021} investigated tissue growth driven by a diffusing morphogen, where domain evolution results from an elasticity-based growth mechanism. More recently, Abels, Blawid and Dolzmann ~\cite{Abels2026} proposed a morphoelastic model for stress-driven growth in the presence of nutrients, combining a reaction--diffusion description of nutrient transport with the multiplicative decomposition of the deformation gradient into elastic and growth contributions. These studies highlight the role of chemical species in regulating growth and mechanical deformation, and illustrate different mathematical frameworks for coupling transport processes, elasticity, and domain evolution in biological tissues.

\medskip

A more comprehensive description of growing biological tissues can be obtained by combining reaction–diffusion or transport processes with morphoelastic growth, so that biochemical fields regulate local growth while the resulting deformation and stresses may, in turn, affect transport and growth kinetics. Such coupled formulations have been considered in several works, including \cite{AmbrosiMollica2002,GarikipatiEtAl2004,AmbrosiGuillou2007,XueEtAl2016,CollisEtAl2017,WalkerEtAl2023,DavoliNikStefanelli2023}, among others. These studies provide different frameworks for coupling nutrient or chemical transport, mechanical deformation, and the evolution of the tissue reference configuration.

\medskip

\medskip
\noindent\textbf{Mechanical modelling of the arterial wall.}
The arterial wall is a heterogeneous, anisotropic, and nearly incompressible
soft tissue whose mechanical response is governed by the nonlinear interaction
between elastin, collagen fibres, and smooth-muscle cells.  The comprehensive
treatment of vascular solid mechanics by Humphrey~\cite{humphrey2002} remains
the standard reference for the constitutive description of arterial tissue,
covering hyperelastic frameworks, residual stresses, and the coupling between
wall mechanics and intraluminal pressure.  The constrained mixture model of
Humphrey and Rajagopal~\cite{humphrey2002b} provides a rigorous
continuum-mechanical basis for describing growth and remodelling of arterial
tissue as the net result of the turnover of individual constituents, each
governed by its own deposition stretch and degradation rate.  In the
pathological setting, the progressive stiffening of the arterial wall due to
lipid accumulation, calcification, and fibrous-cap formation alters local stress
distributions in a manner that may promote plaque rupture, as documented in
detail by Stary et al.~\cite{stary1995} and analyzed mechanically by Holzapfel,
Sommer, and Regitnig~\cite{holzapfel2004} using layer-specific constitutive
models.  The interaction between arterial compliance and pulsatile blood flow,
as described in Nichols and O'Rourke~\cite{nichols2011}, gives rise to
fluid--structure interaction effects whose mathematical treatment requires the
coupling of Navier--Stokes equations in the lumen with elasticity equations in
the wall; see, e.g., Quarteroni, Tuveri, and Veneziani~\cite{quarteroni2000}
and the review by Formaggia, Quarteroni, and Veneziani~\cite{formaggia2009}.
While the present work does not account for full fluid--structure interaction,
it captures the essential mechanical response of the wall through a linearized
elasticity framework whose parameters can be related to the effective Lam\'e
coefficients of the composite arterial tissue.

\medskip
\noindent\textbf{Coupled elasticity--diffusion systems.}
The coupling between mechanical deformation and diffusive transport in elastic
solids arises in a variety of physical and biological contexts, including
thermoelasticity, poroelasticity, and the swelling of hydrogels.  In the
context of biological tissues, such couplings are essential for describing how
concentration gradients induce internal stresses and how mechanical deformation
modifies transport properties.  The generalized thermo-diffusion framework of
Hobiny and Abbas~\cite{HobinyAbbas2022}, from which the elasticity--diffusion system
of the present paper is derived, provides a rigorous continuum-mechanical
formulation in which displacement, concentration, and chemical potential are
governed by a coupled system of PDEs obtained from thermodynamic consistency
arguments.  The classical Biot theory of poroelasticity~\cite{biot1941}
provides an analogous framework for fluid-saturated porous media, and its
extensions to biological tissues --- including articular cartilage and
interstitial tissue --- are reviewed in Mow et al.~\cite{mow1980} and Simon~\cite{simon1992}.  More recently, mechanochemical models have been proposed in
which reaction--diffusion systems are posed on deformable domains and the domain
evolution is driven by the mechanical response to concentration fields; see
Mascheroni et al.~\cite{mascheroni2018} and the references therein.  The present
paper contributes to this line of research by formulating a two-dimensional
model specific to atherosclerotic plaque growth, in which the concentration of
macrophages drives wall deformation through a coupling term in the Lam\'e
system, and by establishing an exact mass--volume relation that characterizes
the global deformation of the domain in terms of boundary displacements.

\medskip

\section{Coupled Elasticity--Diffusion Model in 2D}
We consider a two-dimensional coupled system describing the interaction between plane-strain elasticity and diffusion. The model is inspired by thermo-diffusion formulations in elasticity in paper \cite{HobinyAbbas2022},
where mechanical deformation is coupled with both temperature and concentration.
Such elasticity--diffusion couplings appear in several physical and biological applications where concentration gradients induce deformation and mechanical stress affects the diffusion process \cite{nowacki1976,weitsman1987,sharma2011}.
In the present work, we simplify this framework by neglecting thermal effects and
consider only the coupling between elasticity and diffusion.

The unknowns of the system are the displacement vector $u = (u_x, u_y)$, the concentration $C(x,y,t)$ of the diffusing species, the chemical potential $p(x,y,t)$, and the stress tensor $\sigma$.

The governing equations are given by
\begin{align}
\mu \Delta u_x + (\lambda + \mu) \left( \frac{\partial^2 u_x}{\partial x^2} + \frac{\partial^2 u_y}{\partial x \partial y} \right)
- \gamma_2 \frac{\partial C}{\partial x} &= \rho \frac{\partial^2 u_x}{\partial t^2}, \label{eq:2.1} \\[6pt]
\mu \Delta u_y + (\lambda + \mu) \left( \frac{\partial^2 u_y}{\partial y^2} + \frac{\partial^2 u_x}{\partial x \partial y} \right)
- \gamma_2 \frac{\partial C}{\partial y} &= \rho \frac{\partial^2 u_y}{\partial t^2}, \label{eq:2.2} \\[6pt]
D_0 \Delta p + F(C) &= \frac{\partial C}{\partial t}.\label{eq:2.3}
\end{align}

Equations \eqref{eq:2.1}--\eqref{eq:2.2}  describe the balance of linear momentum in the elastic medium.
The terms involving the Lamé coefficients $\lambda$ and $\mu$ correspond to the
classical elastic response of the material. The additional term involving the gradient
of the concentration introduces a coupling between mechanics and diffusion, meaning
that spatial variations of $C$ generate internal forces and induce deformation.
The right-hand side represents inertial effects, where $\rho$ denotes the density.

Equation \eqref{eq:2.3} describes the diffusion process. The term $D_0 \Delta p$ represents
diffusion driven by the gradient of the chemical potential $p$, while the function $F(C)$
models possible reaction or source terms. The right-hand side corresponds to the time
variation of the concentration.

The stress tensor is defined by
\begin{align}
\sigma_{xx} &= 2\mu \frac{\partial u_x}{\partial x} + \lambda \left(\frac{\partial u_x}{\partial x} + \frac{\partial u_y}{\partial y}\right) - \beta_2 C, \\
\sigma_{xy} &= \sigma_{yx} = \mu \left( \frac{\partial u_x}{\partial y} + \frac{\partial u_y}{\partial x} \right), \\
\sigma_{yy} &= 2\mu \frac{\partial u_y}{\partial y} + \lambda \left(\frac{\partial u_x}{\partial x} + \frac{\partial u_y}{\partial y}\right) - \beta_2 C.
\end{align}

These relations correspond to a linear elastic constitutive law, where the stress depends
on the strain through the Lamé coefficients. The additional term involving $C$ reflects
the influence of the concentration on the mechanical response of the medium.

Finally, the chemical potential is given by
\begin{equation}
p = -\beta_2 \left( \frac{\partial u_x}{\partial x} + \frac{\partial u_y}{\partial y} \right) + bC.
\end{equation}

This relation shows that the chemical potential depends on both the concentration and the deformation of the medium.
The term involving the divergence of $u$ represents the effect of volume changes on the chemical potential. As a result, the mechanical deformation influences the diffusion process.

\subsection{Numerical results for domain deformation under a fixed concentration distribution}

\begin{figure}[H]
\centering

\begin{subfigure}{0.3\textwidth}
\centering
\includegraphics[width=\linewidth, trim={30cm 15cm 0cm 0cm}, clip]{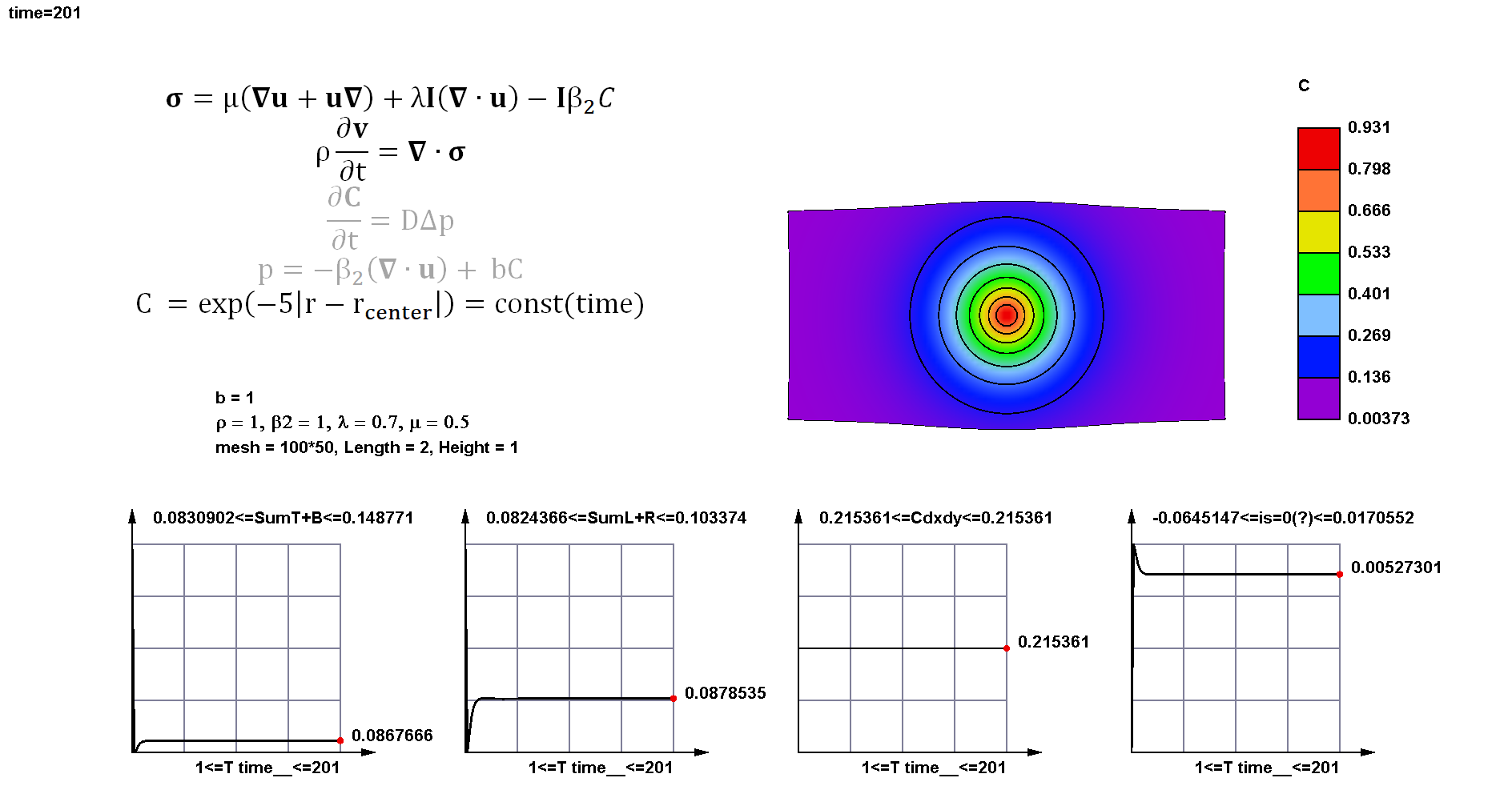}
\caption{$C$ in $2D$ at $t=201$}
\end{subfigure}
\hfill
\begin{subfigure}{0.3\textwidth}
\centering
\includegraphics[width=\linewidth, trim={30cm 15cm 0cm 0cm}, clip]{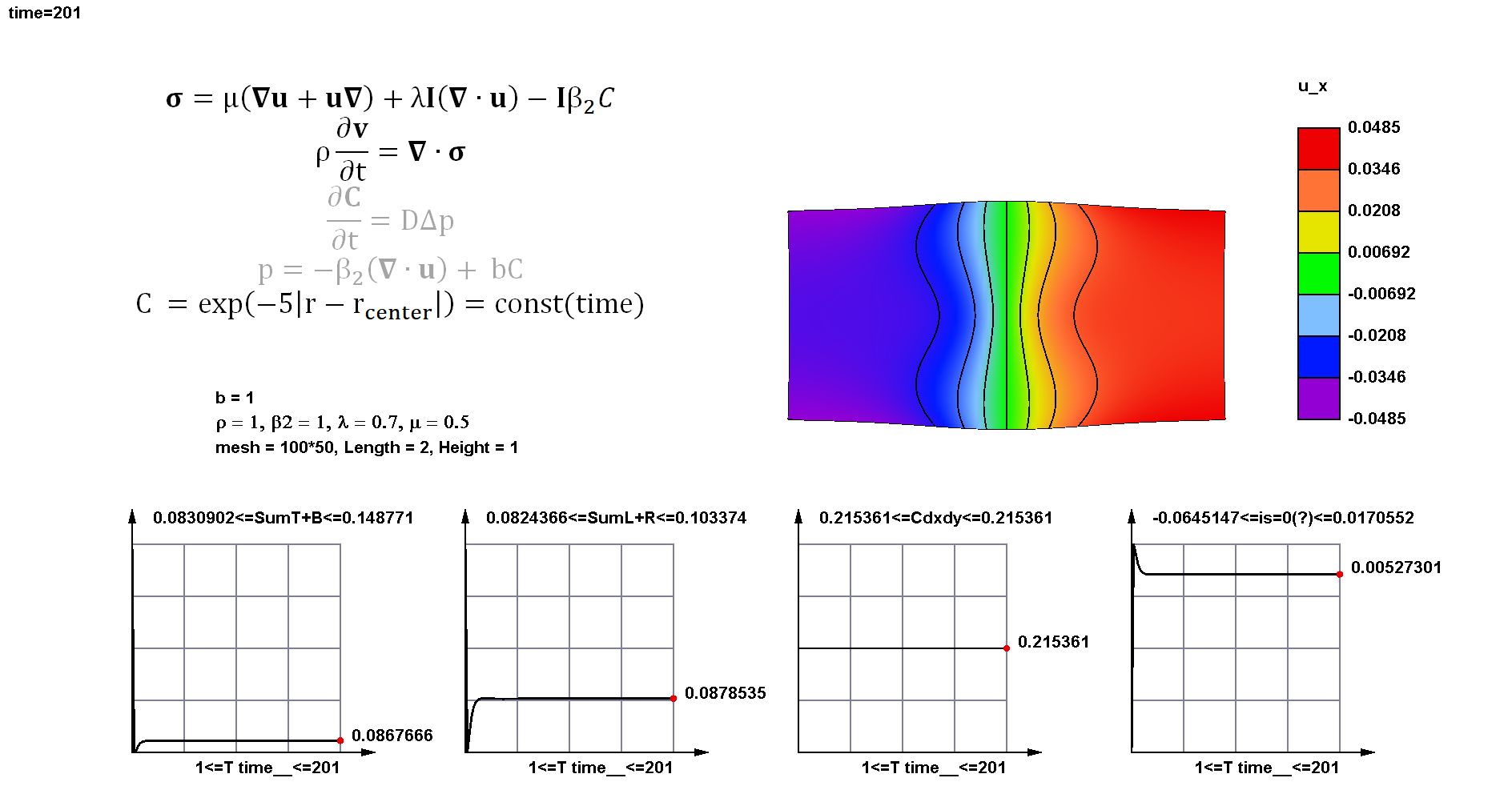}
\caption{$u_x$ in $2D$ at $t=201$}
\end{subfigure}
\hfill
\begin{subfigure}{0.3\textwidth}
\centering
\includegraphics[width=\linewidth, trim={30cm 15cm 0cm 0cm}, clip]{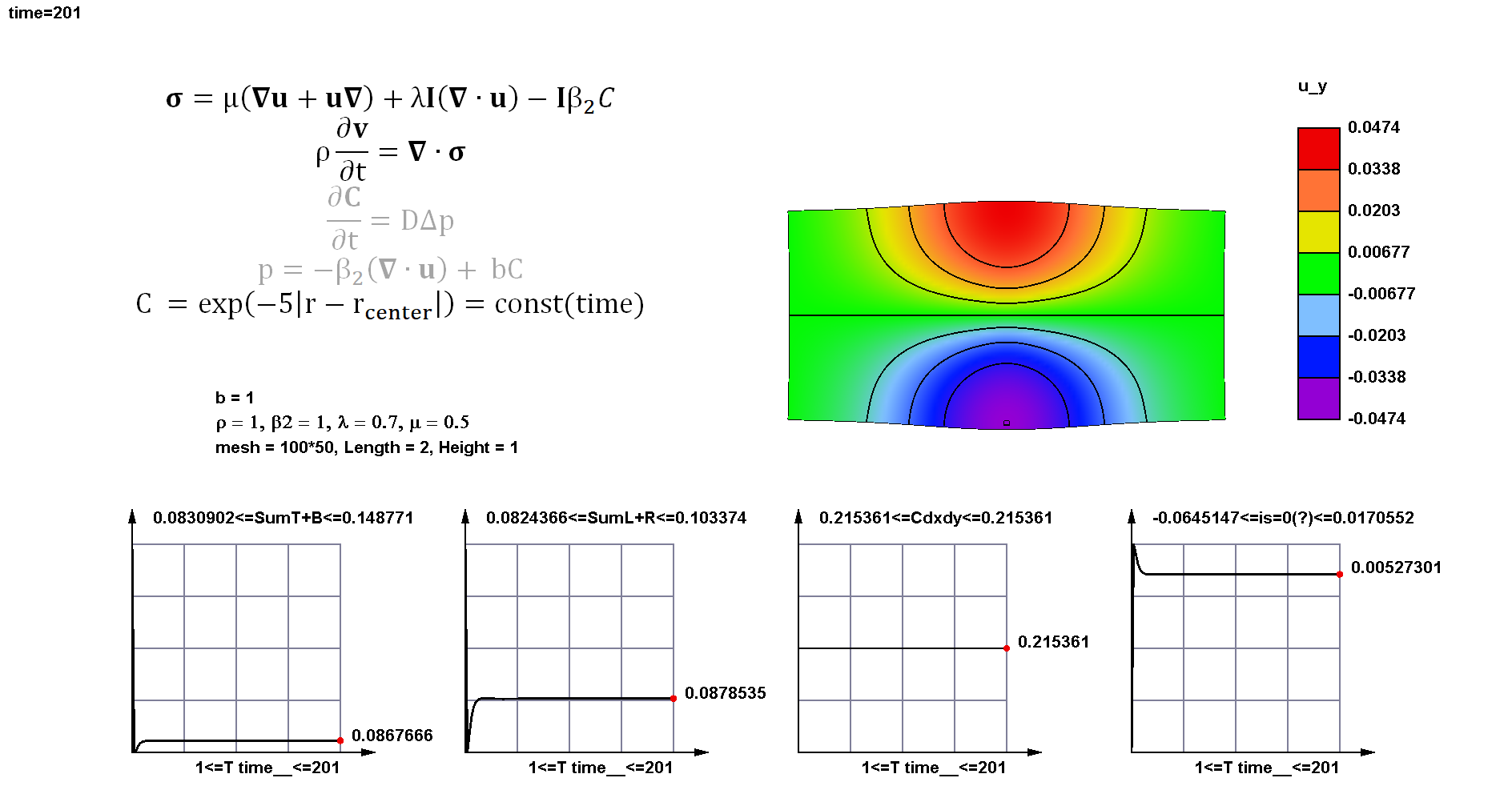}
\caption{$u_y$ in $2D$ at $t=201$}
\end{subfigure}

\vspace{0.5cm}

\begin{subfigure}{0.3\textwidth}
\centering
\includegraphics[width=\linewidth, trim={30cm 15cm 0cm 0cm}, clip]{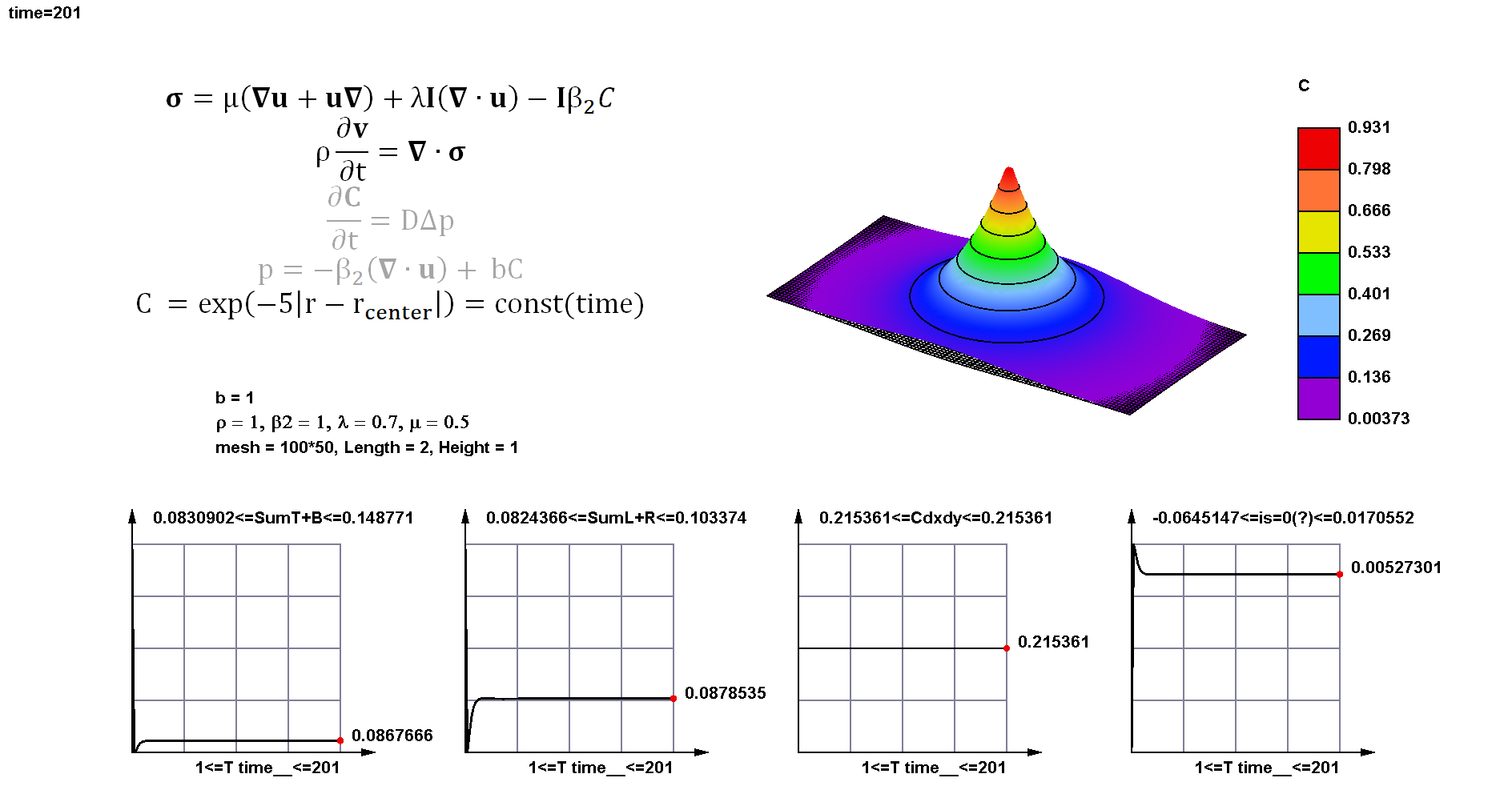}
\caption{$C$ in $3D$ at $t=201$}
\end{subfigure}
\hfill
\begin{subfigure}{0.3\textwidth}
\centering
\includegraphics[width=\linewidth, trim={30cm 15cm 0cm 0cm}, clip]{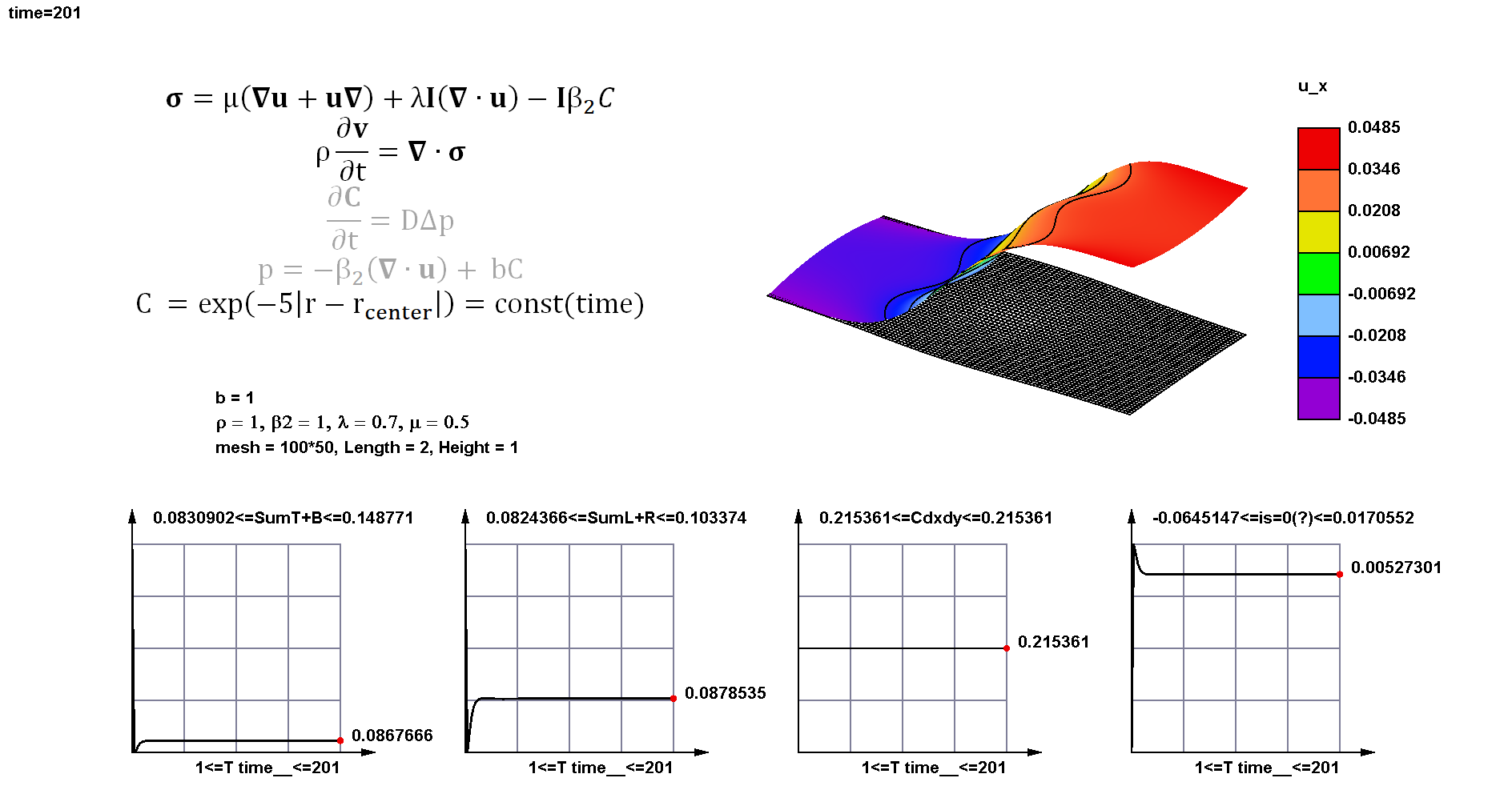}
\caption{$u_x$ in $3D$ at $t=201$}
\end{subfigure}
\hfill
\begin{subfigure}{0.3\textwidth}
\centering
\includegraphics[width=\linewidth, trim={30cm 15cm 0cm 0cm}, clip]{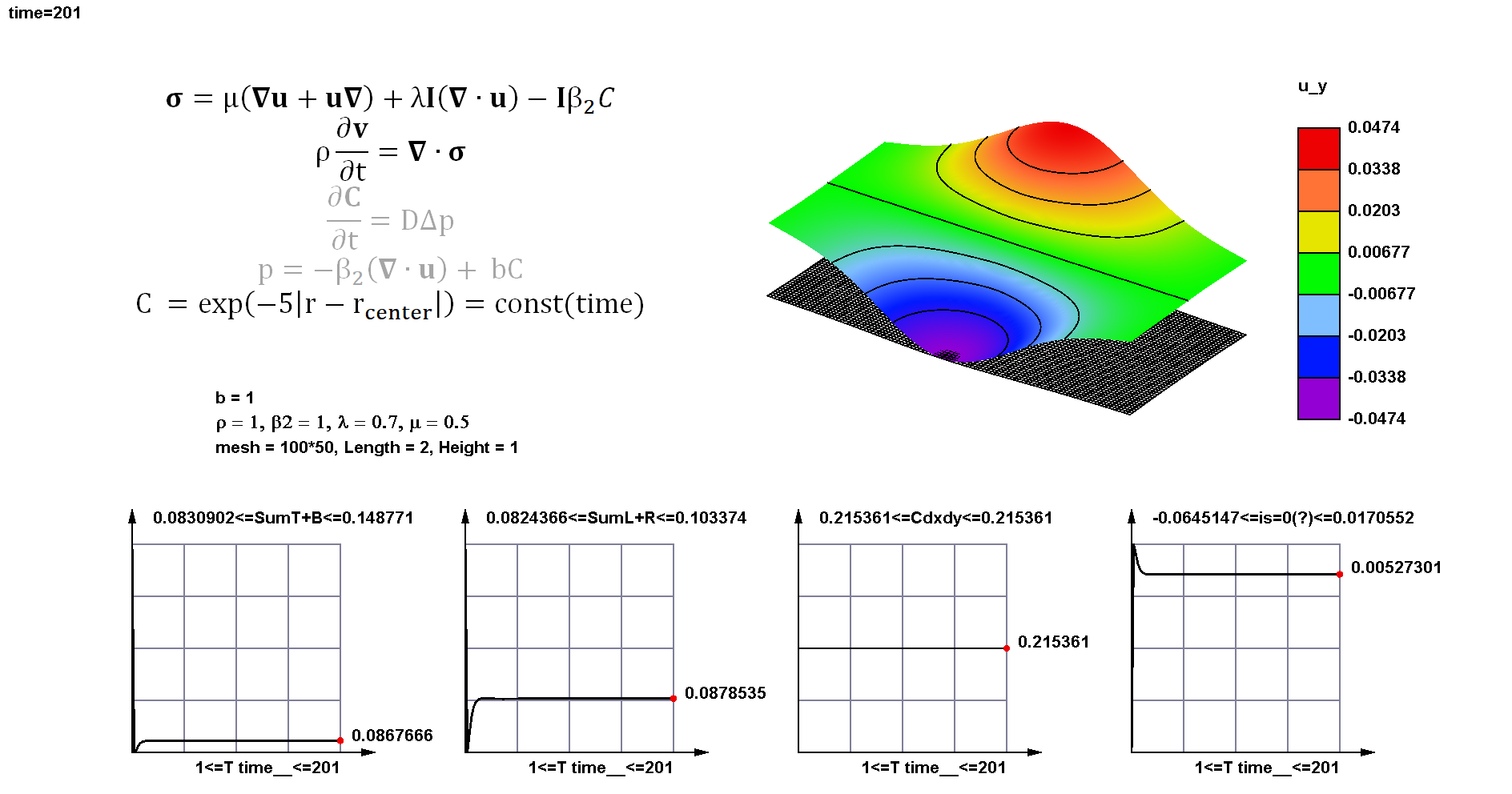}
\caption{$u_y$ in $3D$ at $t=201$}
\end{subfigure}

\caption{Numerical results under zero-stress boundary conditions and a prescribed concentration distribution fixed in time. Panels (a) and (d) show the concentration $C$, panels (b) and (e) show the horizontal displacement $u_x$, and panels (c) and (f) show the vertical displacement $u_y$, in two- and three-dimensional representations, respectively. The deformed configuration illustrates the mechanical response of the domain to the prescribed spatial concentration distribution.}
\label{fig:beta1}
\end{figure}

Figure \ref{fig:beta1} shows the numerical results for the implementation of system \eqref{eq:2.1}--\eqref{eq:2.3} with zero-stress boundary conditions. In this test case, the concentration distribution is prescribed and kept fixed in time in order to illustrate how the domain deforms in response to a given spatial concentration profile. The initial domain is a rectangular domain of height $H=1$ and length $L=2$. The parameters are
\[
\mu = 0.5, \qquad \lambda = 0.7, \qquad \rho = 1, \qquad \gamma_2=\beta=1.
\]
The concentration $C$ is constant in time and is given by
\[
C(x,y)=\exp\bigl(-5|r-r_{\text{center}}|\bigr).
\]

The results in Figure \ref{fig:beta1} show that the prescribed nonuniform concentration distribution generates a nonuniform displacement field and consequently a deformation of the domain. The horizontal displacement $u_x$ varies across the domain in response to the spatial variation of $C$, while the vertical displacement $u_y$ produces a displacement of the upper and lower boundaries. The resulting deformed configuration reflects the spatial distribution of the prescribed concentration. Since the parameter $\beta$ controls the strength of the mechanical response to the concentration field, decreasing $\beta$ reduces the magnitude of the induced displacements and, consequently, the deformation of the domain. These results illustrate the direct effect of the concentration field on the mechanical deformation of the domain, even though the concentration itself is kept fixed in time.

\subsection{Quasi-stationary formulation in 2D}

To simplify the analysis of the coupled system, we consider a quasi-stationary approximation
for the mechanical part. This corresponds to assuming that the mechanical response is
much faster than the diffusion process. Accordingly, in equations (\ref{eq:2.1})-(\ref{eq:2.2}) the inertial terms
$\frac{\partial^2 u_x}{\partial t^2}$ and $\frac{\partial^2 u_y}{\partial t^2}$ can be
neglected, and the medium is assumed to instantaneously reach mechanical equilibrium.

Under this assumption, the balance of linear momentum reduces to
\begin{align}
\gamma_2 \frac{\partial C}{\partial x} &= \mu \Delta u_x + (\lambda+\mu)\left( \frac{\partial^2 u_x}{\partial x^2} + \frac{\partial^2 u_y}{\partial x \partial y}\right), \label{eq:2.7} \\
\gamma_2 \frac{\partial C}{\partial y} &= \mu \Delta u_y + (\lambda+\mu)\left( \frac{\partial^2 u_y}{\partial y^2} + \frac{\partial^2 u_x}{\partial x \partial y}\right). \label{eq:2.8}
\end{align}

Differentiating equation \eqref{eq:2.7} with respect to $x$, equation \eqref{eq:2.8} with
respect to $y$, and summing the resulting expressions, we obtain
\begin{equation}
\gamma_2 \Delta C = (\lambda+2\mu)\left( \frac{\partial}{\partial x}\Delta u_x + \frac{\partial}{\partial y}\Delta u_y \right). \label{eq:2.9}
\end{equation}

Using the relation for the chemical potential, we deduce
\begin{equation}
\Delta p = -\beta_2 \left( \frac{\partial}{\partial x}\Delta u_x + \frac{\partial}{\partial y}\Delta u_y \right) + b\Delta C = \gamma \Delta C, \label{eq:2.10}
\end{equation}
where
\[
\gamma = b -\frac{ \beta_2 \gamma_2}{\lambda+2\mu}.
\]

This relation allows us to express the chemical potential directly in terms of the
concentration and eliminates it from the system. As a result, the coupled model can be
written in the reduced form
\begin{align}
\mu \Delta u_x + (\lambda + \mu) \left( \frac{\partial^2 u_x}{\partial x^2} + \frac{\partial^2 u_y}{\partial x \partial y} \right) - \gamma_2 \frac{\partial C}{\partial x} &= 0, \label{eq:1.11} \\
\mu \Delta u_y + (\lambda + \mu) \left( \frac{\partial^2 u_y}{\partial y^2} + \frac{\partial^2 u_x}{\partial x \partial y} \right) - \gamma_2 \frac{\partial C}{\partial y} &= 0. \label{eq:1.12}
\end{align}

After this reduction, the evolution of the concentration is governed by a reaction--diffusion
equation of the form
\begin{equation}
\frac{\partial C}{\partial t} = D \Delta C + F(C), \label{eq:1.13}
\end{equation}
where
\[
D=D_0\gamma
=D_0\frac{b-\beta_2\gamma_2}{\lambda+2\mu}
\]
is the effective diffusion coefficient, and $F(C)$ represents reaction
or source terms.

Equation \eqref{eq:1.13} is of reaction-diffusion type. After the reduction, this equation
is coupled with the elasticity system \eqref{eq:1.11}--\eqref{eq:1.12}, so that the
concentration influences the deformation of the medium, while the deformation, in turn,
affects the diffusion process.

Under the quasi-stationary assumption, the system reduces to:
\begin{equation}
\left\{
\begin{aligned}
\mu \Delta u_x + (\lambda+\mu)\left(\frac{\partial^2 u_x}{\partial x^2}
+ \frac{\partial^2 u_y}{\partial x \partial y}\right)
- \gamma_2 \frac{\partial C}{\partial x} = 0, \\[6pt]
\mu \Delta u_y + (\lambda+\mu)\left(\frac{\partial^2 u_y}{\partial y^2}
+ \frac{\partial^2 u_x}{\partial x \partial y}\right)
- \gamma_2 \frac{\partial C}{\partial y} = 0, \\[6pt]
\frac{\partial C}{\partial t} = D \Delta C + F(C), \\[6pt]
\sigma_{xx} = 2\mu \frac{\partial u_x}{\partial x} + \lambda \left( \frac{\partial u_x}{\partial x} + \frac{\partial u_y}{\partial y} \right) - \beta_2 C, \\
\sigma_{xy} = \sigma_{yx} = \mu \left( \frac{\partial u_x}{\partial y} + \frac{\partial u_y}{\partial x} \right), \\
\sigma_{yy} = 2\mu \frac{\partial u_y}{\partial y} + \lambda \left( \frac{\partial u_x}{\partial x} + \frac{\partial u_y}{\partial y} \right) - \beta_2 C.
\end{aligned}
\right.
\label{eq:reducedsystem}
\end{equation}

We consider a two-dimensional domain with a moving upper boundary defined by
\[
D(t) = \left\{ (x, y) \mid 0 \le x \le L, \ 0 \le y \le H(x, t) \right\},
\]
where $H(x,t)$ describes the evolution of the upper boundary over time. This allows us
to account for the deformation of the domain induced by the interaction between
diffusion and elasticity.

The boundary conditions are given as follows.

At $x=0$ and $x=L$,
\begin{equation}
   \frac{\partial C}{\partial x} = 0, \quad u_x = 0, \quad \frac{\partial u_y}{\partial x} = 0,
\end{equation}
which correspond to zero flux of concentration through the lateral boundaries.
The horizontal displacement is constrained, while vertical sliding is allowed.

At $y=0$,
\begin{equation}
\frac{\partial C}{\partial y} = 0, \quad u_x = 0, \quad u_y = 0,
\end{equation}
which represents a no-flux condition for the concentration and a fixed boundary
for the displacement.

At $y = H(x,t)$,
\begin{equation}
\frac{\partial C}{\partial n} = f(C), \quad
\sigma_{xx} n_x + \sigma_{yx} n_y = 0, \quad
\sigma_{xy} n_x + \sigma_{yy} n_y = 0. \end{equation}
Here $n=(n_x,n_y)$ denotes the outward unit normal vector. The condition on $C$
models an influx of cells or chemical species through the boundary, while the
mechanical conditions correspond to the absence of external forces on the upper boundary.

Overall, the bottom boundary is fixed, the lateral boundaries allow vertical motion,
and the upper boundary is free of normal stresses.
\subsubsection{Numerical approaches}
In this section, we consider two approaches for the numerical simulation of the coupled system, based on how the deformation of the domain is handled.

In the first approach, the computational domain is kept fixed, with $H$ taken as a constant.
The equations are solved on the reference rectangular domain. The deformation is not incorporated
into the geometry itself, but is instead evaluated through the displacement field. In other words,
the change in the domain at each time step is measured relative to the original domain, without
updating the computational mesh.

In the second approach, the computational domain evolves in time. More precisely, the upper boundary
is described by the function $H=H(x,t)$, which is updated at each time step according to the displacement.
This leads to a moving domain and the equations are solved on this evolving geometry. In this case, both the diffusion process and
the mechanical deformation are computed on the deformed domain, which provides a more accurate
description of the coupling between concentration and domain deformation.

The fixed-domain approach provides a simpler approximation that avoids the complexity of a moving boundary,
while the moving-domain approach describes the full geometric evolution of the system.

\subsubsection{Numerical results for the quasi-stationary system}
In this section, we present numerical results obtained using the two approaches described above, namely the fixed-domain approach and the moving-domain approach. The purpose is to compare their behavior and illustrate the effect of domain deformation on the evolution of the concentration and displacement.

In all simulations, we consider the parameters
\[
D = 0.1, \quad \mu = 0.5, \quad \lambda = 0.7, \quad \gamma_2 = 1, \quad \beta_2 = 1, \quad L= 2, \quad H(x,0) = 1.
\]
The initial concentration is taken as
\[
C_0(x,y)=\frac{1}{100} \exp \Big( -\big((x - 1)^2 + (y - 0.5)^2\big) \Big),
\]
and the boundary function is given by $f(C)=0.2C$. The time step is chosen as $dt = 0.001$, and the computational mesh consists of $23385$ elements.

\begin{figure}[h!]
\centering

\begin{subfigure}{0.45\textwidth}
    \centering
    \includegraphics[width=\linewidth]{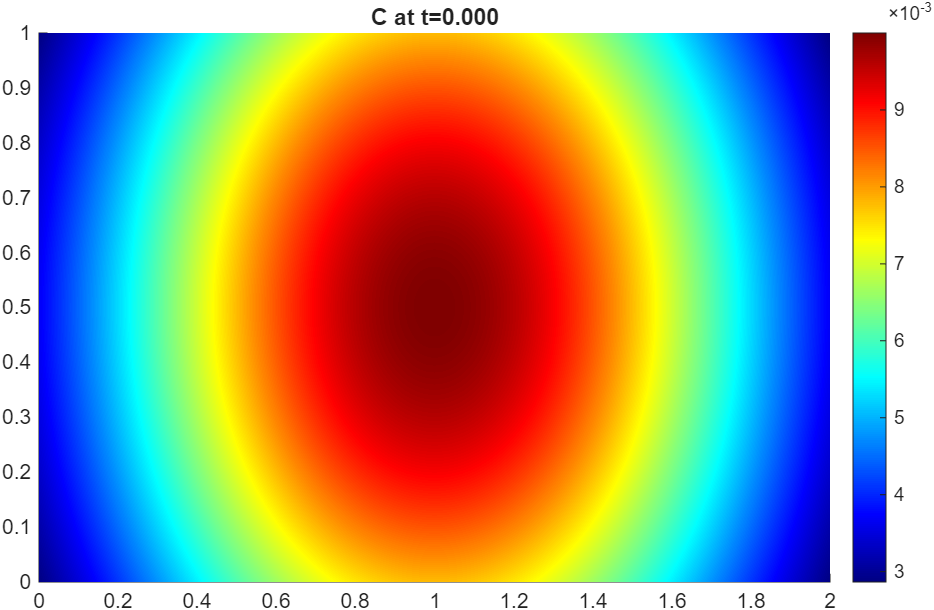}
    \caption{$C$ at $t=0$}
\end{subfigure}
\hfill
\begin{subfigure}{0.45\textwidth}
    \centering
    \includegraphics[width=\linewidth]{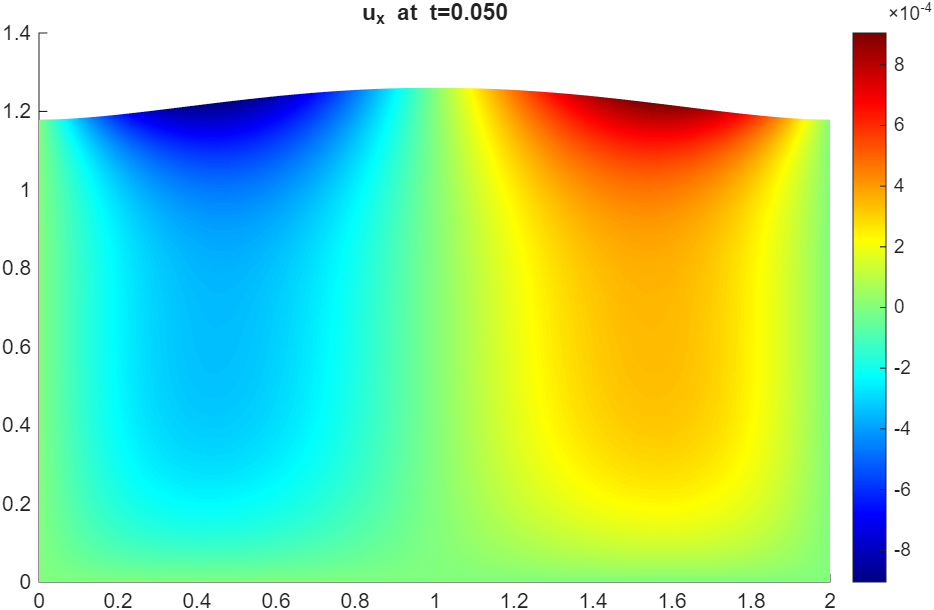}
    \caption{$u_x$ at $t=0.05$}
\end{subfigure}

\vspace{0.5cm}

\begin{subfigure}{0.45\textwidth}
    \centering
    \includegraphics[width=\linewidth]{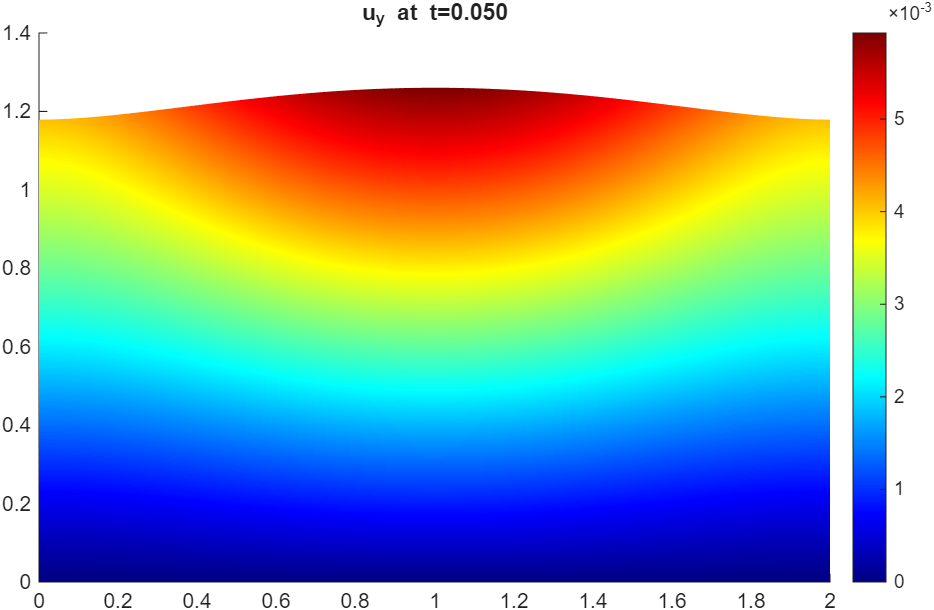}
    \caption{$u_y$ at $t=0.05$}
\end{subfigure}
\hfill
\begin{subfigure}{0.45\textwidth}
    \centering
    \includegraphics[width=\linewidth]{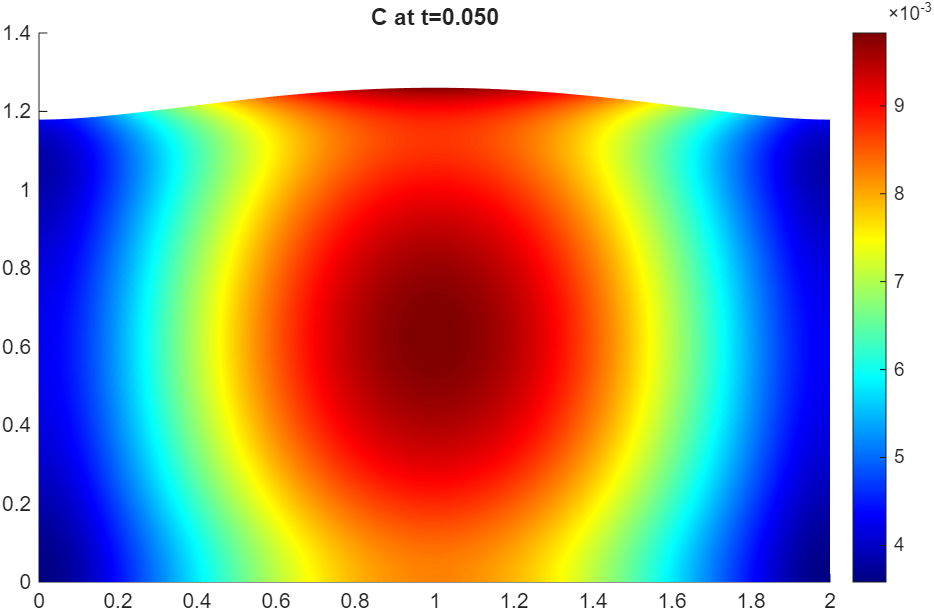}
    \caption{$C$ at $t=0.05$}
\end{subfigure}
\caption{Results obtained using the moving-domain approach. Panel (a) shows the initial concentration distribution $C$ at $t=0$. Panel (b) shows the horizontal displacement $u_x$ at $t=0.05$, with opposite signs on the two sides of the domain, indicating horizontal deformation in opposite directions. Panel (c) shows the vertical displacement $u_y$ at $t=0.05$, which is positive and more pronounced near the upper part of the domain, resulting in an upward deformation of the upper boundary. Panel (d) shows the concentration distribution $C$ at $t=0.05$ on the resulting deformed domain.
}
\label{M2D}
\end{figure}

The results of the simulations obtained using the moving-domain approach are presented in Figure \ref{M2D}. The concentration distribution generates both horizontal and vertical displacements. The horizontal displacement $u_x$ has opposite signs on the two sides of the domain, corresponding to an outward deformation, while the vertical displacement $u_y$ increases toward the upper boundary and produces an upward displacement of this boundary. Consequently, the initially rectangular domain becomes deformed, with a more pronounced displacement near the central part of the upper boundary, where the concentration is higher.

Table 1 reports the area change at selected time steps for the fixed domain and moving domain approaches. For small deformations, the two approaches give close results. As the deformation increases, however, the difference between the two approaches becomes more pronounced. This is due to the fact that, in the fixed-domain approach, the equations are always solved on the original reference domain, whereas in the moving-domain approach the geometry is updated at each time step. Therefore, the effect of the evolving geometry becomes increasingly significant as the deformation grows.

\begin{table}[h!]
\centering
\begin{tabular}{c|c|c}
\hline
$t$ & Area change with fixed domain approach & Area change with moving domain approach \\
\hline
0.001 & 0.008343 & 0.008343 \\
0.005 & 0.008171 & 0.008323 \\
0.010 & 0.008162 & 0.008505 \\
0.015 & 0.008162 & 0.008702 \\
0.020 & 0.008166 & 0.008907 \\
0.025 & 0.008171 & 0.009119 \\
0.030 & 0.008177 & 0.009337 \\
0.035 & 0.008184 & 0.009562 \\
0.040 & 0.008191 & 0.009793 \\
0.045 & 0.008198 & 0.010031 \\
0.050 & 0.008206 & 0.010276 \\
\hline
\end{tabular}
\caption{Comparison of the area change for the fixed domain and moving domain approaches at selected time values.}
\label{comparison}
\end{table}

\subsection{Mass--Volume Relation}\label{sec:massvolume}

In this section, we establish a relation between the total mass of the concentration field and the deformation of the domain. Such a relation is relevant from a physical point of view, as it directly connects the amount of diffusing species in the medium to the resulting domain deformation, and it also provides a useful criterion for assessing the accuracy of the numerical simulations. More precisely, we show that, under the assumptions of linear elasticity and small deformations, the total mass can be expressed in terms of the volume change of the domain. The derivation is based on integrating the elasticity equations over the domain and using the stress-free boundary conditions. The calculations are presented in Appendix~\ref{appendix:massvolume}. Here, we summarize the main steps and state the resulting relation.

Consider the elasticity system \eqref{eq:2.7}--\eqref{eq:2.8} in the rectangular domain
\[
0<x<L,\qquad 0<y<H,
\]
with displacement field $(u_x,u_y)=(u,v)$ and $\gamma_2=\beta_2=\beta$. We assume the stress-free boundary conditions
\[
\sigma_{xx}=0 \quad \text{on } x=0,\qquad
\sigma_{yy}=0 \quad \text{on } y=0,
\]
together with
\[
\sigma_{xy}=0
\quad\text{on}\quad
x=0,\;x=L,\;y=0,\;\text{and}\;y=H.
\]

Integrating the elasticity equations over the domain and applying the above boundary conditions, we obtain
\begin{equation}
\beta \int_0^H \int_0^L C(x,y)\,dx\,dy
=
(2\mu+\lambda)\int_0^H (u(L,y)-u(0,y))\,dy
+\lambda\int_0^L (v(x,H)-v(x,0))\,dx,
\label{eq:mass1}
\end{equation}
and
\begin{equation}
\beta \int_0^L \int_0^H C(x,y)\,dx\,dy
=
\lambda\int_0^H (u(L,y)-u(0,y))\,dy
+
(\lambda+2\mu)\int_0^L (v(x,H)-v(x,0))\,dx.
\label{eq:mass2}
\end{equation}

The detailed derivation of these identities is presented in Appendix~\ref{appendix:massvolume}. Comparing \eqref{eq:mass11} and \eqref{eq:mass2} yields
\begin{equation}
\int_0^H (u(L,y)-u(0,y))\,dy
=
\int_0^L (v(x,H)-v(x,0))\,dx,
\label{eq:equalboundary}
\end{equation}
showing that the total boundary displacement in the horizontal direction equals that in the vertical direction.

For small deformations, the variation of the domain volume is approximated by the first-order contribution of the boundary displacement. The geometric derivation is given in Appendix~\ref{appendix:massvolume}, where it is shown that the quadratic terms are negligible. Consequently, the volume change is given by
\begin{equation}
\Delta V
=
\int_0^H (u(L,y)-u(0,y))\,dy
+
\int_0^L (v(x,H)-v(x,0))\,dx.
\label{eq:DeltaV}
\end{equation}

Using \eqref{eq:equalboundary}, we obtain
\begin{equation}
\Delta V
=
2\int_0^L (v(x,H)-v(x,0))\,dx.
\label{eq:DeltaV2}
\end{equation}

Substituting this expression into either \eqref{eq:mass1} or \eqref{eq:mass2} gives the following mass--volume relation:
\begin{equation}
\Delta V
=
\frac{\beta}{\lambda+\mu}
\int_0^L\int_0^H C(x,y)\,dx\,dy.
\label{linearc}
\end{equation}

Therefore, under the assumptions of linear elasticity and small deformations, the total volume change is proportional to the total mass of the concentration field, with proportionality constant
\[
\frac{\beta}{\lambda+\mu}.
\]
This relation will also be used below to assess the consistency of the numerical results.

The rectangular derivation presented here can be extended to an arbitrary
domain. As shown in Appendix~\ref{app:general_mass_volume}, the resulting
mass--volume relation is independent of the domain geometry when the entire
boundary is traction-free. For the mixed boundary conditions used in the
simulations, an additional contribution arises from the reaction traction on
the clamped boundary.

\subsubsection{Numerical results}
In this section, we consider system \eqref{eq:reducedsystem} together with the boundary conditions used in the derivation presented in Appendix A in order to verify numerically the mass--volume relation. These boundary conditions are given by:

\begin{equation}
\begin{cases}
x=0: & \dfrac{\partial C}{\partial n} = 0, \quad \sigma_{xx} = 0, \quad \sigma_{xy} = 0, \\[6pt]
y=0: & \dfrac{\partial C}{\partial n} = 0, \quad \sigma_{yy} = 0, \quad \sigma_{xy} = 0, \\[6pt]
x=L: & \dfrac{\partial C}{\partial n} = 0, \quad \sigma_{xy} = 0, \\[6pt]
y=H: & \dfrac{\partial C}{\partial n} = f(C), \quad \sigma_{xy} = 0.
\end{cases}
\label{BCcalculations}
\end{equation}
The goal of this section is to examine whether the linear correlation between the mass and the volume in equation (\ref{linearc}) is obtained numerically. The values of the parameters are as follows: $\mu =40000$, $\lambda =60000$, $\gamma_2 =4$,
$\beta_2=4$, $L=5$ and $H=6$. The mesh consists of $152815$  elements and the time step is $dt=0.01$.  To analyze this behavior, we consider three different initial conditions:
\[
C_1(x,y) =
\cos\!\left(\frac{\pi(x - 3)}{6}\right)
\sin\!\left(\frac{3\pi y}{5}\right),
\]

\[
C_2(x,y) =
10\,\exp\!\left(
- \left[
(x - 3)^2 + (y - 2.5)^2
\right]
\right),
\]
and

\[
C_3(x,y) =
\cos\!\left(\frac{3\pi(x - 3)}{6}\right)
\sin\!\left(\frac{\pi y}{5}\right).
\]

\begin{figure}[H]
\centering

\begin{subfigure}[b]{0.32\textwidth}
    \centering
    \includegraphics[width=\textwidth]{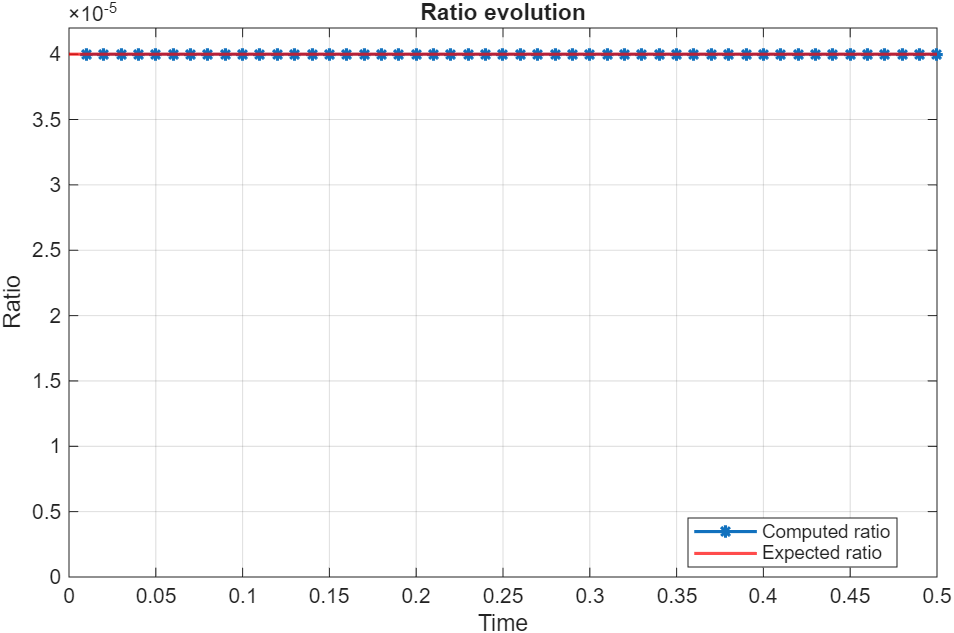}
    \caption{$C_1$ as initial concentration}
\end{subfigure}
\hfill
\begin{subfigure}[b]{0.32\textwidth}
    \centering
    \includegraphics[width=\textwidth]{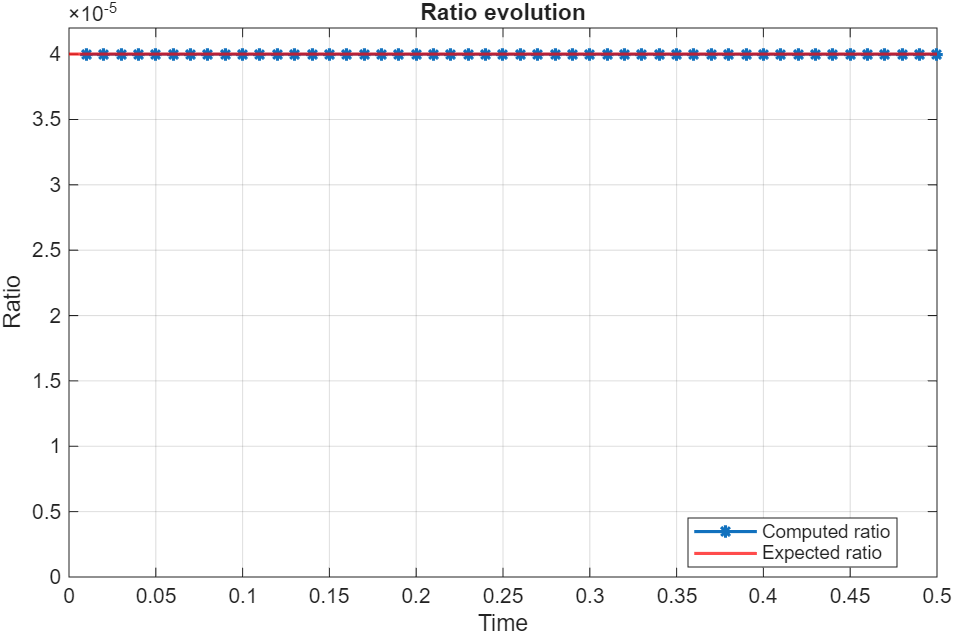}
    \caption{$C_2$ as initial concentration}
\end{subfigure}
\hfill
\begin{subfigure}[b]{0.32\textwidth}
    \centering
    \includegraphics[width=\textwidth]{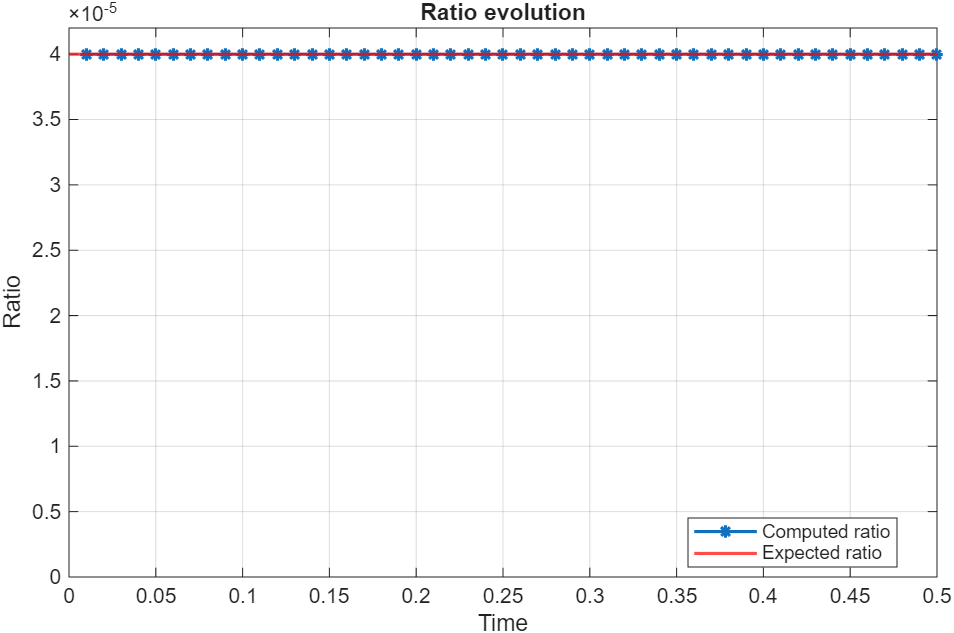}
    \caption{$C_3$ as initial concentration}
\end{subfigure}

\caption{Numerical simulations of \eqref{eq:reducedsystem} with boundary conditions \eqref{BCcalculations} for $D=0$, $f(C)=0$, and $F(C)=0$.}
\label{fig:Case1}
\end{figure}

\begin{figure}[H]
\centering

\begin{subfigure}[b]{0.32\textwidth}
    \centering
    \includegraphics[width=\textwidth]{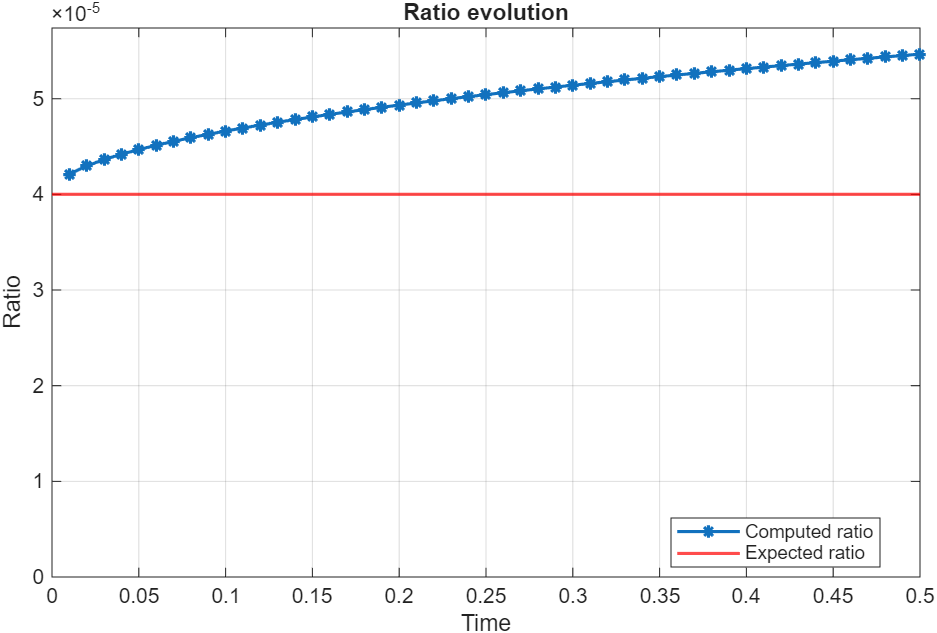}
    \caption{$C_1$ as initial concentration}
\end{subfigure}
\hfill
\begin{subfigure}[b]{0.32\textwidth}
    \centering
    \includegraphics[width=\textwidth]{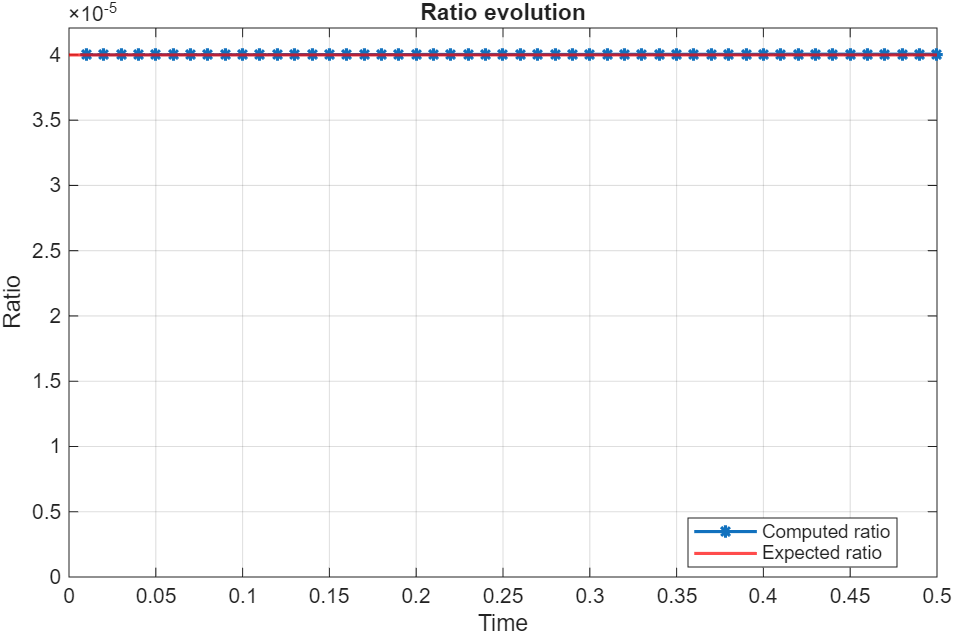}
    \caption{$C_2$ as initial concentration}
\end{subfigure}
\hfill
\begin{subfigure}[b]{0.32\textwidth}
    \centering
    \includegraphics[width=\textwidth]{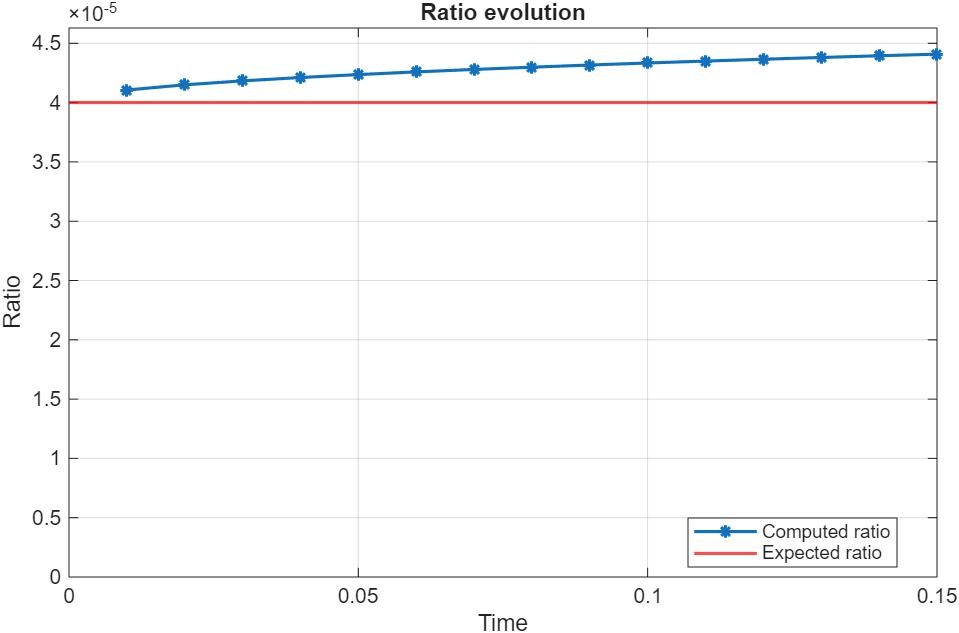}
    \caption{$C_3$ as initial concentration}
\end{subfigure}

\caption{Numerical simulations of \eqref{eq:reducedsystem} with boundary conditions \eqref{BCcalculations} for $D=10^{-2}$, $f(C)=0$, and $F(C)=0$.}
\label{fig:Case2}
\end{figure}

\begin{figure}[H]
\centering

\begin{subfigure}[b]{0.32\textwidth}
    \centering
    \includegraphics[width=\textwidth]{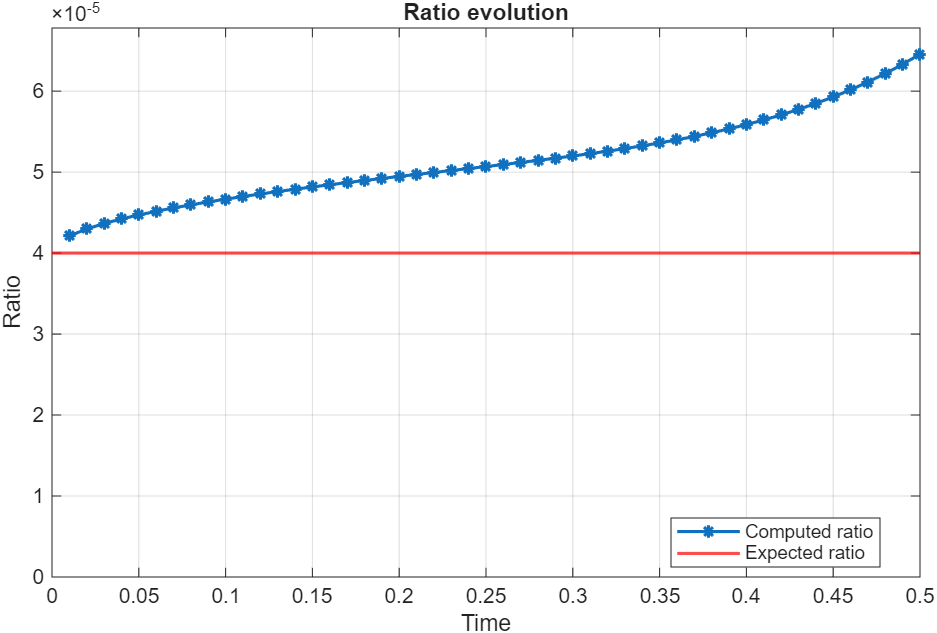}
    \caption{$C_1$ as initial concentration}
\end{subfigure}
\hfill
\begin{subfigure}[b]{0.32\textwidth}
    \centering
    \includegraphics[width=\textwidth]{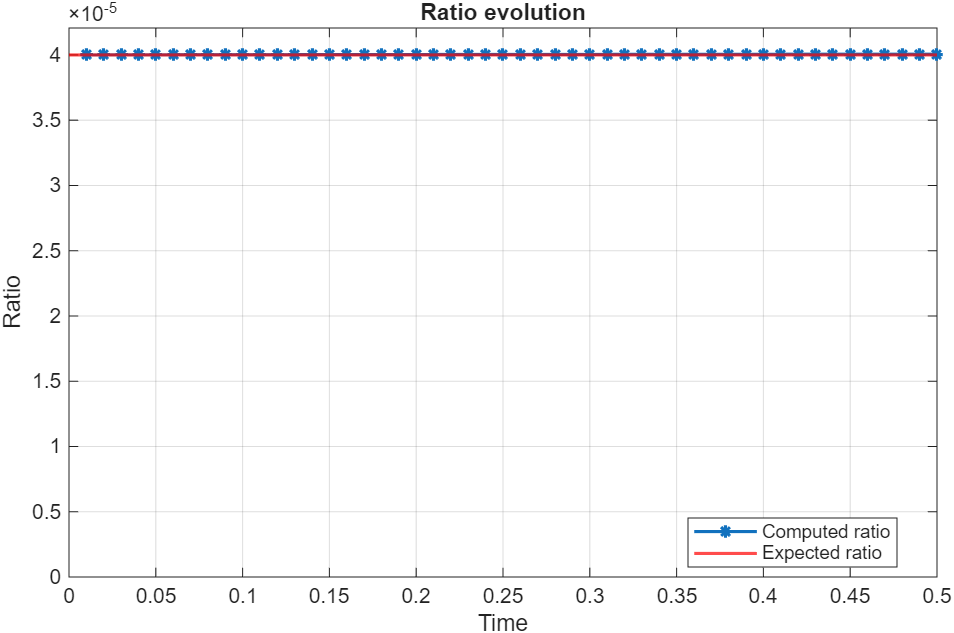}
    \caption{$C_2$ as initial concentration}
\end{subfigure}
\hfill
\begin{subfigure}[b]{0.32\textwidth}
    \centering
    \includegraphics[width=\textwidth]{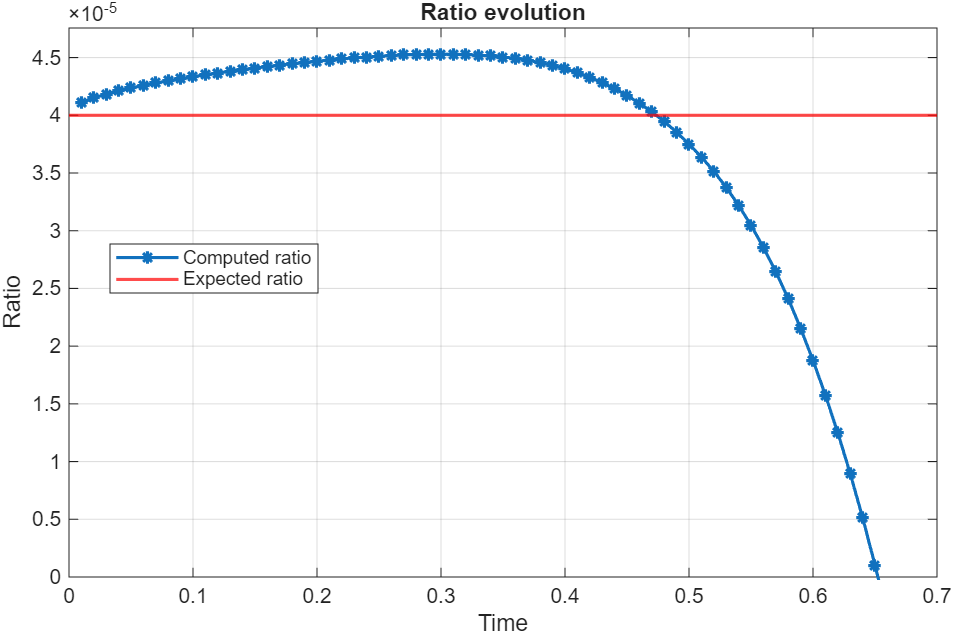}
    \caption{$C_3$ as initial concentration}
\end{subfigure}

\caption{Numerical simulations of \eqref{eq:reducedsystem} with boundary conditions \eqref{BCcalculations} for $D=10^{-2}$, $f(C)=0.5C$, and $F(C)=0$.}
\label{fig:Case3}
\end{figure}

\begin{figure}[H]
\centering

\begin{subfigure}[b]{0.32\textwidth}
    \centering
    \includegraphics[width=\textwidth]{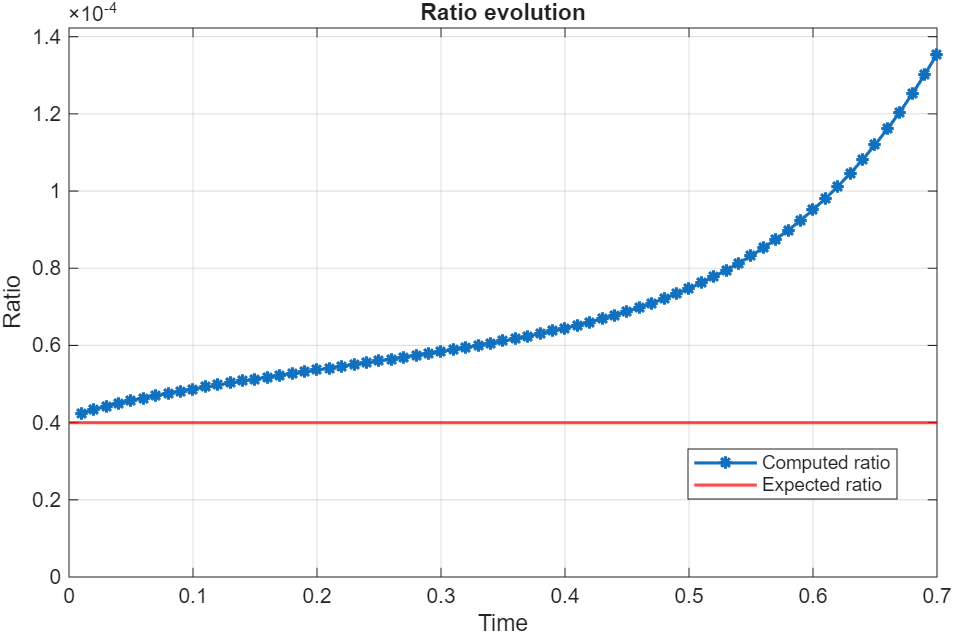}
    \caption{$C_1$ as initial concentration}
\end{subfigure}
\hfill
\begin{subfigure}[b]{0.32\textwidth}
    \centering
    \includegraphics[width=\textwidth]{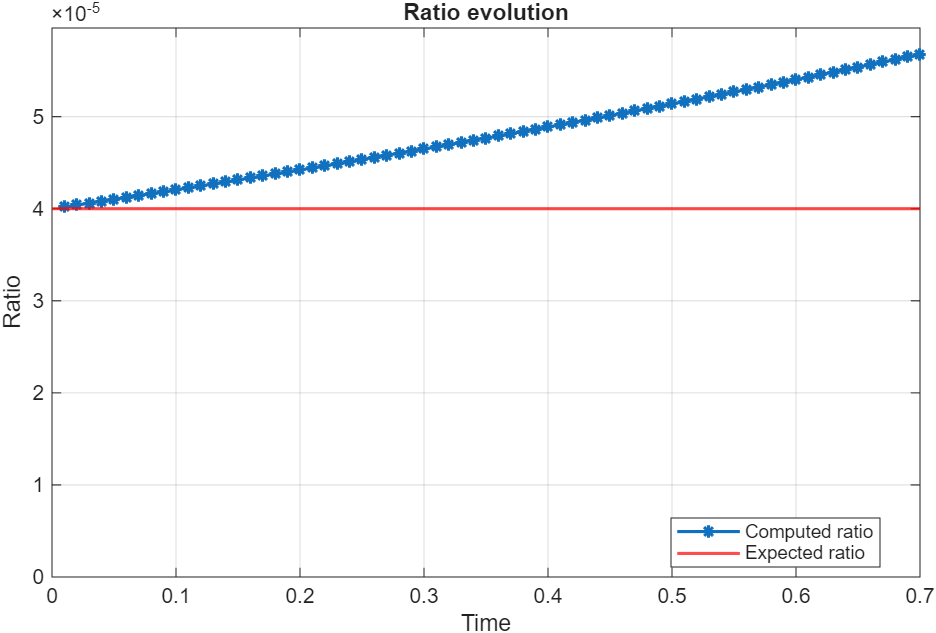}
    \caption{$C_2$ as initial concentration}
\end{subfigure}
\hfill
\begin{subfigure}[b]{0.32\textwidth}
    \centering
    \includegraphics[width=\textwidth]{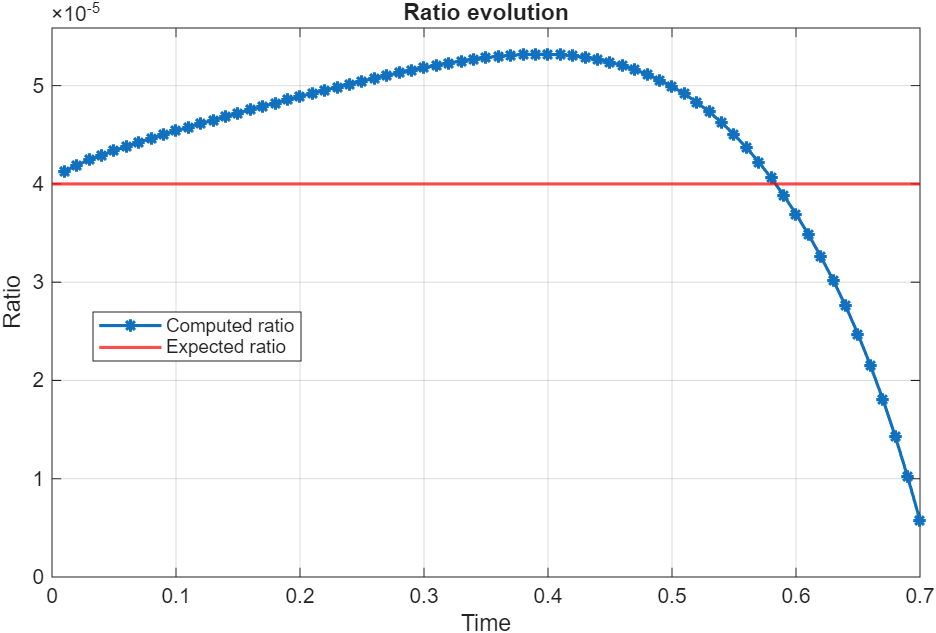}
    \caption{$C_3$ as initial concentration}
\end{subfigure}

\caption{Numerical simulations of \eqref{eq:reducedsystem} with boundary conditions \eqref{BCcalculations} for $D=10^{-2}$, $f(C)=0.5C$, and $F(C)=0.5C$.}
\label{fig:Case4}
\end{figure}

The ratio computed in the numerical simulations corresponds to the area
deformation divided by the double integral of the concentration $C$ over
the domain, namely,
\begin{equation}
\frac{\Delta V}
{\displaystyle\int_0^L\int_0^H C(x,y),dx,dy}.
\label{eq:computed_ratio}
\end{equation}
In Figures \ref{fig:Case1}--\ref{fig:Case4}, this quantity is referred to as the computed ratio and is compared with the expected ratio, denoted by $r$, which is obtained from the analytical mass--volume relation and is given by
\begin{equation}
r=\frac{\beta}{\lambda+\mu}.
\label{eq:expected_ratio}
\end{equation}

Thus, the agreement between the computed ratio and the expected value $r$
provides a numerical verification of the analytical mass--volume relation.

We first consider the case where $D=0$, $f(C)=0$, and $F(C)=0$, corresponding to the absence of diffusion, reaction or source terms, and boundary fluxes. In this case, the concentration equation reduces to
$\displaystyle\frac{\partial C}{\partial t}=0,$ and therefore the concentration does not depend on time, i.e., $C(x,y,t)=C(x,y,0).$ Thus, the initial concentration distribution remains unchanged throughout the simulation. System (\ref{eq:reducedsystem}), together with the boundary conditions (\ref{BCcalculations}), is implemented using three different initial concentration distributions. The numerical simulations show that the linear correlation remains preserved throughout all time steps; see Figure~\ref{fig:Case1}. In addition, the displacement remains very small, with a magnitude of approximately
$4\times10^{-5}$, which corresponds exactly to the value of the ratio $r$ for the
chosen parameter set.

Next, we consider the same system \eqref{eq:reducedsystem} with the same boundary conditions \eqref{BCcalculations}, while introducing diffusion by taking \(D=10^{-2}\), together with \(f(C)=0\) and \(F(C)=0\). The corresponding numerical results are presented in Figure~\ref{fig:Case2}. Since diffusion only redistributes the concentration and no source or boundary flux is present, the integral of \(C\) remains constant. Therefore, the ratio is expected to remain equal to its initial value. The results obtained for the three initial concentrations show that the computed ratios remain very close to the expected value, but increase slightly at each time step. This small variation indicates that the higher-order terms neglected in deriving the linear mass--volume relation are not entirely negligible in this case. Although the displacement remains small, the contribution of these terms relative to the small area change is sufficient to produce the observed deviation.

We then consider the case
$D=10^{-2}$,
$f(C)=0.5C$,
and $F(C)=0$,
thereby introducing a source term in the equation.
The obtained numerical simulations are shown in Figure~\ref{fig:Case3}. This figure shows that the ratio \(r\) is not always preserved. However, a closer examination of the displacement values reveals that as long as the displacement remains of the order of \(10^{-4}\), the ratio obtained stays very close to the expected one. Once the displacement becomes larger (of the order of \(10^{-3}\) or greater), the ratio changes significantly, and the corresponding curve deviates from a straight line.

Finally, we study the case
$D=10^{-2}$,
$f(C)=0.5C$,
and $F(C)=0.5C$,
which includes both a source term and a boundary flux contribution.
The corresponding results are displayed in Figure~\ref{fig:Case4}. The simulations show similar behaviors to those obtained in Figure~\ref{fig:Case3}.

These numerical results show that the linear correlation obtained in
\eqref{linearc} is preserved only as long as the deformation remains small.
Indeed, the derivation of \eqref{linearc} relies on a linear elasticity framework
and on the approximation of the volume change by its first-order contribution.
In the calculation of the deformed area, higher-order terms involving products of
the displacement components and their derivatives were neglected. This is justified
when the displacement is very small, since these nonlinear terms are then negligible
compared with the leading-order terms.

However, when the displacement increases, the deformation of the boundary is no
longer accurately described by the linear approximation. In this case, the actual change in area contains nonlinear geometric contributions that are not included in the theoretical derivation of the ratio. As a result, the computed ratio defined in \eqref{eq:computed_ratio} is no longer expected to remain equal to the expected ratio $r$ given by \eqref{eq:expected_ratio}.

Therefore, the observed deviation of the numerical ratio from the expected value
does not contradict the analytical calculation. Rather, it indicates that the
linear mass-volume relation is valid in the small-deformation regime, while larger
displacements introduce nonlinear elastic and geometric effects that modify the
area deformation and consequently change the computed ratio.

\section{Coupled Elasticity--Diffusion Model in 1D}
The two-dimensional model provides a general framework for describing the coupling between diffusion and elasticity in a deformable medium. We now introduce a one-dimensional reduction of this model, which can be viewed as a particular case of the two-dimensional system. Despite its simpler formulation, it retains the essential coupling between concentration and mechanical deformation.

In this section, we study a one-dimensional elastic body coupled with a diffusion process. The model describes the interaction between the mechanical deformation of the medium and the distribution of a chemical concentration along a one-dimensional interval.

We consider the problem in the interval $[0,L]$, where the unknowns depend on the spatial variable $x$ and time $t$. The model involves the displacement $u(x,t)$, the concentration $C(x,t)$, the stress $\sigma(x,t)$, and the chemical potential $p(x,t)$. After presenting the complete system, we neglect relaxation effects and assume mechanical equilibrium in order to derive a simplified one-dimensional model. We then discuss some qualitative properties of the solutions, including the deformation of the elastic body and the motion of the free boundary induced by concentration variations.

\subsection{Formulation}

The influence of temperature and concentration distributions on the elastic body has been extensively studied in the literature \cite{nowacki1976,weitsman1987,sharma2011,HobinyAbbas2022} (and references therein). We begin with the following model introduced in \cite{HobinyAbbas2022}:

\begin{equation}
(\lambda + 2\mu)
\frac{\partial^2 u}{\partial x^2}
- \beta_2
\frac{\partial C}{\partial x}
=
\rho
\frac{\partial^2 u}{\partial t^2},
\label{eq:1D-full-mechanics}
\end{equation}

\begin{equation}
D
\frac{\partial^2 p}{\partial x^2}
=
\frac{\partial C}{\partial t}
+
\tau_1
\frac{\partial^2 C}{\partial t^2},
\label{eq:1D-full-diffusion}
\end{equation}

\begin{equation}
\sigma
=
(\lambda + 2\mu)
\frac{\partial u}{\partial x}
-
\beta_2 C,
\label{eq:1D-full-stress}
\end{equation}

\begin{equation}
p
=
-\beta_2
\frac{\partial u}{\partial x}
+
bC.
\label{eq:1D-full-potential}
\end{equation}

Here $u$ denotes the displacement, $C$ the concentration, $\sigma$ the stress, and $p$ the chemical potential. The physical meaning of the coefficients appearing in this model can be found in \cite{HobinyAbbas2022}. This system of equations is considered in the interval $[0,L]$ with some boundary conditions which will be specified below. In what follows, we neglect relaxation terms (second time derivatives) and consider the simplified model:

\begin{equation}
(\lambda + 2\mu)
\frac{\partial^2 u}{\partial x^2}
-
\beta_2
\frac{\partial C}{\partial x}
=
0,
\label{eq:1D-mech-equilibrium}
\end{equation}

\begin{equation}
D
\frac{\partial^2 p}{\partial x^2}
=
\frac{\partial C}{\partial t},
\label{eq:1D-diffusion-potential}
\end{equation}

\begin{equation}
\sigma
=
(\lambda + 2\mu)
\frac{\partial u}{\partial x}
-
\beta_2 C,
\label{eq:1D-stress}
\end{equation}

\begin{equation}
p
=
-\beta_2
\frac{\partial u}{\partial x}
+
bC.
\label{eq:1D-potential}
\end{equation}

Neglecting the time derivative in \eqref{eq:1D-full-mechanics}, we assume that the system is in mechanical equilibrium, that is, no elastic waves are considered.

\subsection{Properties of solutions}

\textbf{Example.} Suppose that $C(x,t)\equiv C_0$. Then from \eqref{eq:1D-mech-equilibrium} we obtain
\[
u(x,t)=k_1(t)x+k_2(t),
\]
and from \eqref{eq:1D-diffusion-potential},
\[
p(x,t)=k_3(t)x+k_4(t),
\]
where $k_i(t)$ are arbitrary functions. Next,
\[
\sigma=(\lambda+2\mu)k_1(t)-\beta_2 C_0,
\]
\[
p=-\beta_2 k_1(t)+bC_0.
\]

Set $u(x,0)=0$ (fixed boundary) and $\sigma(L,t)=0$ (free boundary). Then
\[
k_1(t)=\frac{\beta_2 C_0}{\lambda+2\mu},
\qquad
k_2(t)=0,
\qquad
u(x,t)=\frac{\beta_2 C_0 x}{\lambda+2\mu}.
\]

Hence,
\[
\sigma=0
\]
for all $x$. We obtain that if one boundary is fixed, another one is free and the concentration is constant, then there is no stress in the 1D elastic rode. Its elongation is
\[
u(L,t)=\frac{bC_0L}{\lambda+2\mu}.
\]

\textbf{Motion of free boundary.} Integrating equation \eqref{eq:1D-mech-equilibrium} from $0$ to $x$, we obtain
\[
(\lambda+2\mu)
\frac{\partial u}{\partial x}
-
\beta_2 C(x)
-
\left(
(\lambda+2\mu)
\frac{\partial u}{\partial x}
-
\beta_2 C(x)
\right)_{x=0}
=0.
\]

Assuming that $\sigma(0,t_0)=0$ (free boundary), we conclude from \eqref{eq:1D-stress} that the second term in the previous equality vanishes. Integrating again this equality from $0$ to $L$, we have

\begin{equation}
u(L,t)-u(0,t)
=
\kappa
\int_0^L C(x,t)\,dx
\label{eq:1D-displacement-L}
\end{equation}

where
\[
\kappa=\frac{\beta_2}{\lambda+2\mu}.
\]

Taking into account that
\begin{equation}
L(t)-L_0
=
u(L_0,t)-u(0,t),
\label{eq:1D-length-displacement}
\end{equation}
where $L_0$ is the initial interval length, we obtain the following equality

\begin{equation}
L(t)
=
L_0
+
\kappa
\int_0^{L(t)} C(x,t)\,dx.
\label{eq:1D-free-boundary}
\end{equation}

Relation \eqref{eq:1D-free-boundary} was previously considered in
\cite{abiyounes2024} in the context of a one-dimensional free-boundary model of
atherosclerosis. In that work, this relation was introduced empirically, without being derived from a mechanical model, whereas in the present work it is obtained as a consequence of the elasticity equations. In \cite{abiyounes2024}, the motion of the boundary was determined
by the accumulation of cells within the domain, leading to an integral
relation of the same form as \eqref{eq:1D-free-boundary}. Differentiating
this relation with respect to time gives an evolution law for the free
boundary in which the boundary velocity depends, in general, not only on
the flux through the boundary but also on the production, degradation,
and accumulation of the species inside the domain. Therefore, the resulting
free-boundary problem is, in general, different from a classical Stefan
problem. Under additional simplifying assumptions, however, the evolution law
reduces to a Stefan-type free-boundary condition in which the velocity of
the moving boundary is determined by the flux through the boundary. In
\cite{abiyounes2024}, the moving-boundary problem was transformed to a fixed-domain
problem by a change of variables, and the existence of solutions was
established using a fixed-point argument.

\section{Application to atherosclerosis}
In this section, we model plaque formation in atherosclerosis in a two-dimensional domain.
The model takes into account both the elasticity of the arterial wall and the key players
of the inflammatory process, namely macrophages and cytokines. We consider the coupled
elasticity--diffusion system introduced in \cite{HobinyAbbas2022}, together with the
reaction--diffusion system proposed in \cite{ElKhatib2007} to describe the inflammatory
response in two dimensions. Combining these two approaches, we obtain a system similar to
\eqref{eq:reducedsystem}, in which the equation for the concentration $C$ is replaced by a
system of reaction--diffusion equations describing the inflammatory processes involved
in atherosclerosis.

More precisely, we denote by $u=(u_x,u_y)$ the displacement vector, where $u_x$ and $u_y$
represent the displacement in the $x$ and $y$ directions, respectively. The variable $M$
denotes the concentration of macrophages, $A$ the concentration of cytokines, and $\sigma$
the stress tensor.
We consider the following coupled system:
\begin{equation}
\left\{
\begin{aligned}
\mu \Delta u_x + (\lambda+\mu)\left(\frac{\partial^2 u_x}{\partial x^2}
+ \frac{\partial^2 u_y}{\partial x \partial y}\right)
- \gamma_2 \frac{\partial M}{\partial x} = 0, \\[6pt]
\mu \Delta u_y + (\lambda+\mu)\left(\frac{\partial^2 u_y}{\partial y^2}
+ \frac{\partial^2 u_x}{\partial x \partial y}\right)
- \gamma_2 \frac{\partial M}{\partial y} = 0, \\[6pt]
\frac{\partial M}{\partial t} = d_1 \Delta M - \lambda_1 M, \\[6pt]
\frac{\partial A}{\partial t} = d_2 \Delta A + \frac{\alpha_2 A}{1 + A/\tau_2}\, M - \lambda_2 A,
\\[6pt]
\sigma_{xx} = 2\mu \frac{\partial u_x}{\partial x} + \lambda \left( \frac{\partial u_x}{\partial x} + \frac{\partial u_y}{\partial y} \right) - \beta_2 M, \\
\sigma_{xy} = \sigma_{yx} = \mu \left( \frac{\partial u_x}{\partial y} + \frac{\partial u_y}{\partial x} \right), \\
\sigma_{yy} = 2\mu \frac{\partial u_y}{\partial y} + \lambda \left( \frac{\partial u_x}{\partial x} + \frac{\partial u_y}{\partial y} \right) - \beta_2 M.
\end{aligned}
\right.
\label{eq:reducedsystemathero}
\end{equation}

The first two equations describe the mechanical equilibrium of the arterial wall.
The terms involving the Lamé coefficients $\lambda$ and $\mu$ represent the elastic
response of the material, while the terms involving the gradient of the concentration
account for the effect of macrophages, whose spatial variations induce deformation
of the wall.

The third and fourth equations are reaction--diffusion equations. The diffusion terms
$d_1 \Delta M$ and $d_2 \Delta A$ describe the spatial spread of macrophages and cytokines,
while the terms $-\lambda_1 M$ and $-\lambda_2 A$ represent their degradation. The nonlinear
term
\[
\frac{\alpha_2 A}{1 + A/\tau_2}\, M
\]
models the production of cytokines by macrophages, with a saturation effect.

Finally, the stress tensor $\sigma$ is given by a linear elastic law, where the stress
depends on the strain through the coefficients $\lambda$ and $\mu$. The term
involving $M$ represents the influence of the concentration of macrophages on the mechanical behavior
of the arterial wall.

We consider a two-dimensional domain with a moving upper boundary defined by
\[
D(t) = \left\{ (x, y) \mid 0 \le x \le L, \ 0 \le y \le H(x, t) \right\},
\]
where $H(x,t)$ describes the evolution of the upper boundary over time.

The boundary conditions are given as follows.

At $x=0$ and $x=L$:
\begin{equation}
\frac{\partial M}{\partial x}=0,
\quad
\frac{\partial A}{\partial x}=0, \quad u_x = 0, \quad \frac{\partial u_y}{\partial x} = 0.
\label{RSABC1}
\end{equation}

At $y=0$:
\begin{equation}
   \frac{\partial M}{\partial y}=0,
\quad
\frac{\partial A}{\partial y}=0,\quad u_x = 0, \quad u_y = 0.
\label{RSABC2}
\end{equation}
At $y = H(x,t)$:
\begin{equation}
\frac{\partial M}{\partial y}=\frac{1}{d_1}\left(\frac{\alpha_1 + \beta_1 A}{1 + A/\tau_1}\right),
\quad
\frac{\partial A}{\partial y}=0, \quad
\sigma_{xx} n_x + \sigma_{yx} n_y = 0, \quad
\sigma_{xy} n_x + \sigma_{yy} n_y = 0.
\label{RSABC3}\end{equation}

On the lateral boundaries $x=0$ and $x=L$, and on the bottom boundary $y=0$, homogeneous
Neumann conditions are imposed for $M$ and $A$, meaning that there is no flux of macrophages
or cytokines through these boundaries. The displacement is constrained so that the bottom
boundary is fixed, while the lateral boundaries allow vertical motion.

At the upper boundary, the condition
\[
\frac{\partial M}{\partial y}=\frac{1}{d_1}\left(\frac{\alpha_1 + \beta_1 A}{1 + A/\tau_1}\right)
\]
models the recruitment of monocytes from the blood flow into the intima, which then differentiate into
macrophages. As described in \cite{ElKhatib2007}, this flux depends on the concentration of cytokines at
the boundary: cytokines promote the adhesion and penetration of monocytes through the endothelium.
The term $\alpha_1$ represents a baseline recruitment due to the presence of LDL, while the term
$\beta_1 A$ accounts for the auto-amplification of the recruitment induced by inflammatory cytokines.
The denominator $1 + A/\tau_1$ represents a saturation effect, corresponding to the limitation of
cell recruitment when the concentration becomes large.

The condition $\frac{\partial A}{\partial y}=0$ indicates that there is no flux of cytokines
through the upper boundary, while the stress conditions correspond to the absence of
external forces acting on the boundary.

Without elasticity, the model reduces to a reaction--diffusion system for
the concentration fields on a fixed domain. Such a formulation was studied
in \cite{elkhatib2012}, where a two-dimensional reaction--diffusion model
of atherosclerosis was considered and traveling-wave solutions describing
the propagation of inflammation along the arterial wall were investigated.
This reaction--diffusion model provides a reference case for the present
study, allowing us to investigate how the additional coupling with elasticity
and the resulting mechanical deformation of the domain affect the propagation
of inflammation.

Plaque evolution has also been considered in previous works using
free-boundary formulations of Stefan type. In these models, the motion of
the boundary is associated with the transport, influx, and accumulation of
cells within the plaque. Such a formulation was also considered in
\cite{abiyounes2024}, where the position of the free boundary was related to the
concentrations of cells in the intima. As discussed in Section~4.2, the
resulting free-boundary problem can reduce, under suitable assumptions, to
a Stefan-type problem in which the boundary velocity is determined by the
cell flux. This mechanism is different from the one considered in the
present model where the displacement of the boundary results from the
mechanical deformation of the tissue induced by the coupling between the
concentration fields and elasticity.

\begin{figure}[htbp]
    \centering

    \begin{subfigure}[b]{0.3\textwidth}
        \centering
        \includegraphics[width=\textwidth]{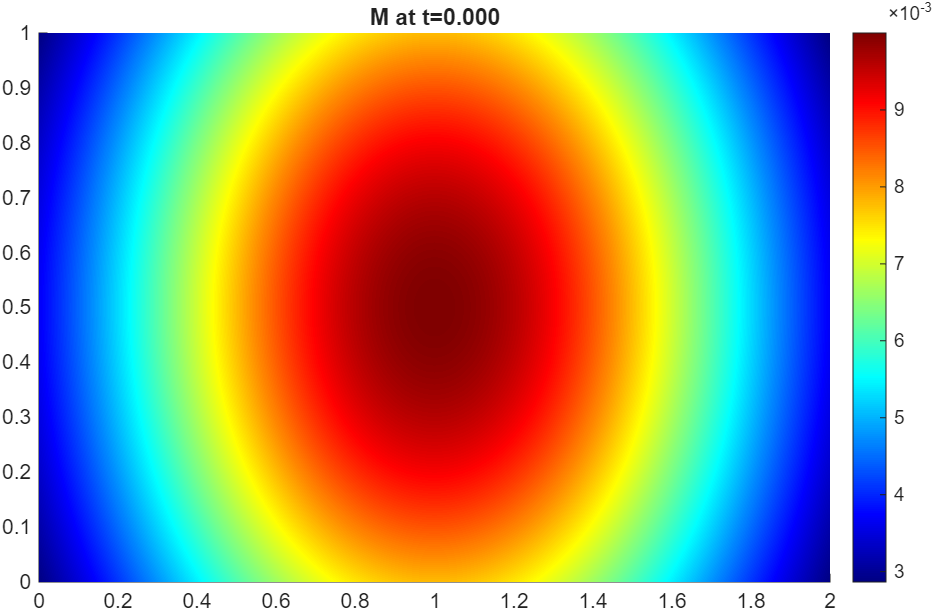}
        \caption{$M$ at $t=0$}
    \end{subfigure}
    \hfill
    \begin{subfigure}[b]{0.3\textwidth}
        \centering
        \includegraphics[width=\textwidth]{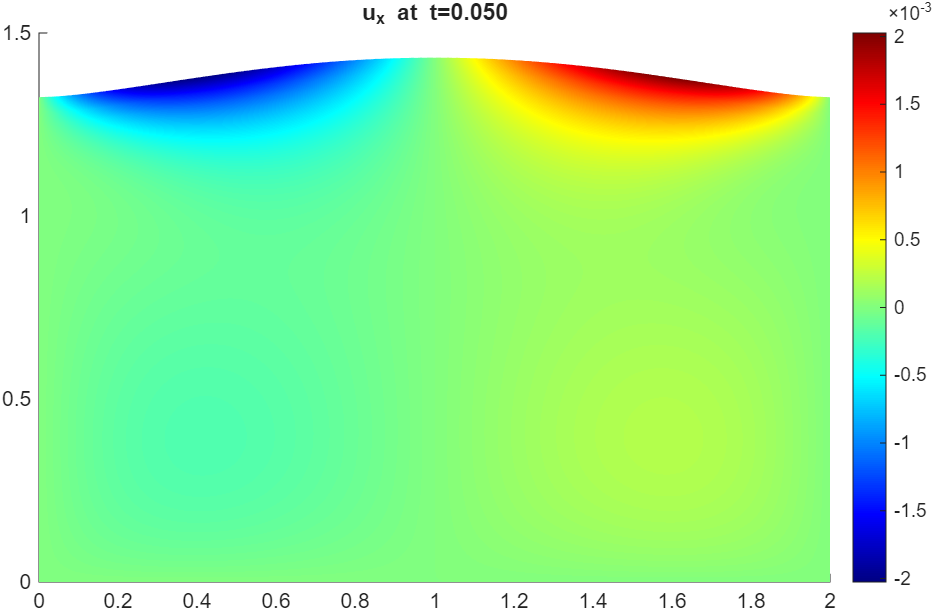}
        \caption{$u_x$ at $t=0.05$}
    \end{subfigure}
    \hfill
    \begin{subfigure}[b]{0.3\textwidth}
        \centering
        \includegraphics[width=\textwidth]{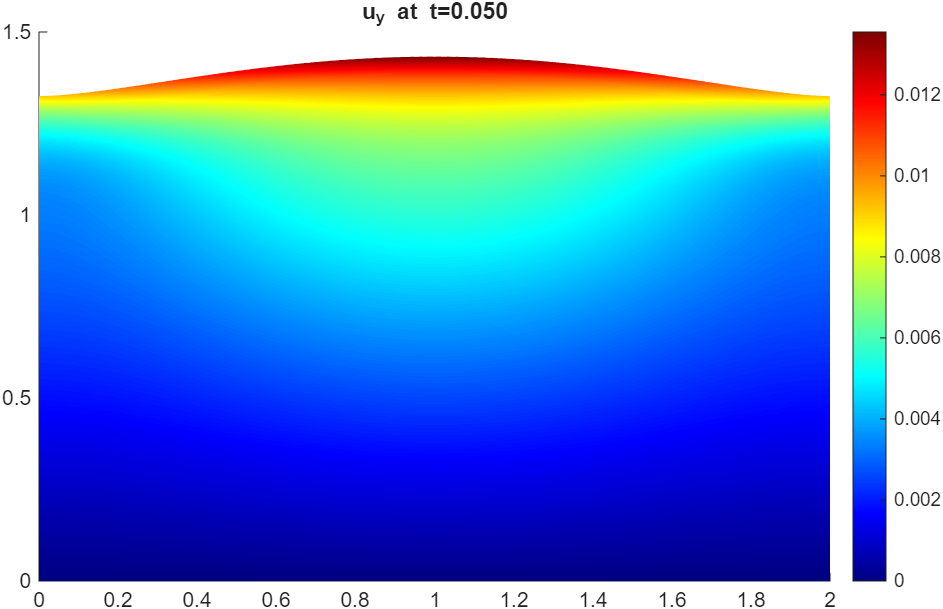}
        \caption{$u_y$ at $t=0.05$}
    \end{subfigure}

    \vspace{0.5cm}

    \begin{subfigure}[b]{0.3\textwidth}
        \centering
        \includegraphics[width=\textwidth]{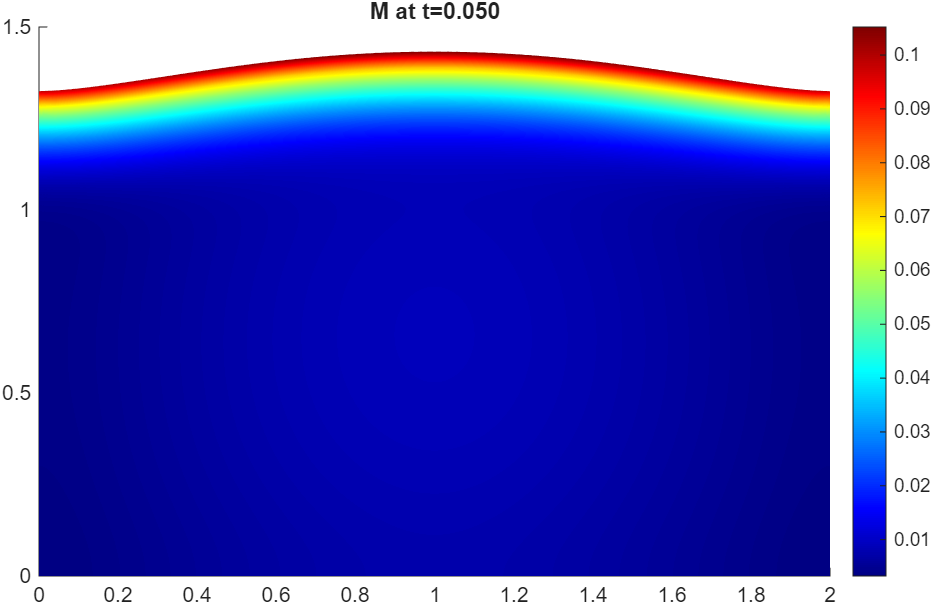}
        \caption{$M$ at $t=0.05$}
    \end{subfigure}
    \hfill
    \begin{subfigure}[b]{0.3\textwidth}
        \centering
        \includegraphics[width=\textwidth]{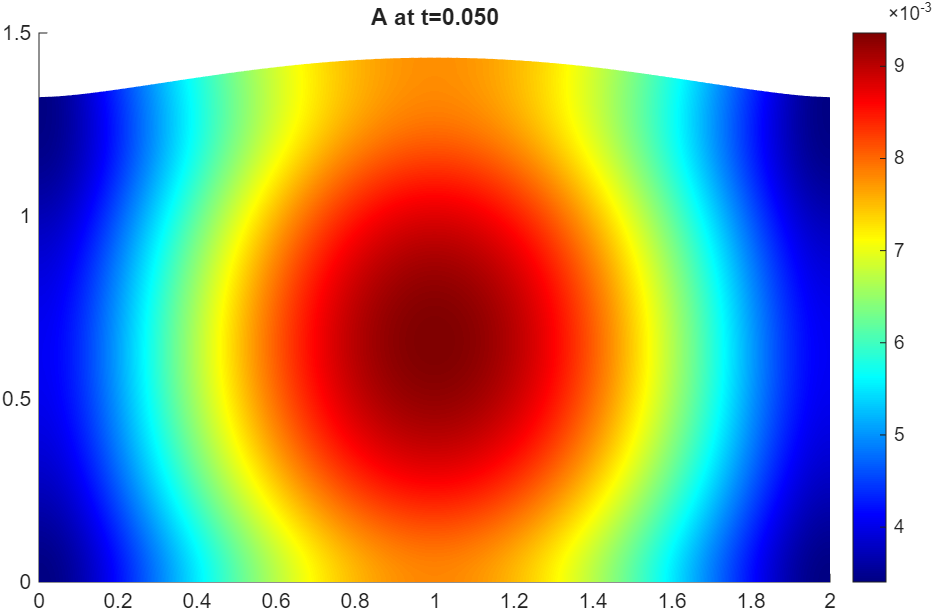}
        \caption{$A$ at $t=0.05$}
    \end{subfigure}
    \hfill
    \begin{subfigure}[b]{0.3\textwidth}
        \centering
        \includegraphics[width=\textwidth]{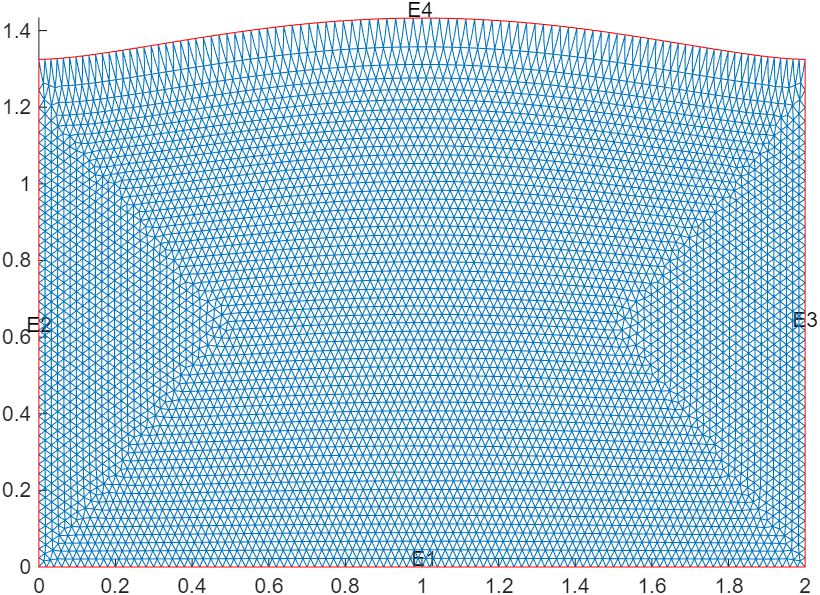}
        \caption{Domain at $t=0.05$}
    \end{subfigure}

    \caption{Numerical solutions in 2D and domain configuration at the indicated times using the moving domain approach.}
    \label{fig:RSAMD2D}
\end{figure}

Figures \ref{fig:RSAMD2D} shows the numerical results obtained by implementing system \eqref{eq:reducedsystemathero} with boundary conditions \eqref{RSABC1}--\eqref{RSABC3} over a rectangular domain of width $2$ and height $1$.

The parameters taken are
$d_1 = 0.1, \quad d_2 = 0.1, \quad
\alpha_1 = 0.1, \quad \alpha_2 = 2,\quad
\beta_1 = 0.1, \quad \lambda_1 = 2, \quad
\lambda_2 = 1, \quad \tau_1 = 1,\quad \tau_2 = 1,\quad
\mu = 0.5, \quad \lambda = 0.7,\quad\gamma_2 = 1, \quad \text{and}\quad \beta_2 = 1.
$ The time step is $\Delta t = 0.001$ and the mesh consists of $23385$ elements. The initial concentrations of $M$ and $A$ are taken as
\[
M_0(x,y) = A_0(x,y) = \frac{1}{100} \exp\left(-\left((x-1)^2 + (y-0.5)^2\right)\right).
\]
We used two approaches for the simulations: the first based on a fixed domain and the second on a moving domain.

Figure~\ref{fig:RSAMD2D} illustrates the response of the coupled system using
the moving-domain approach. Panel (a) shows the initial distribution of
macrophages $M$, which is concentrated around the center of the domain.
At $t=0.05$, the nonuniform distribution of $M$ induces a mechanical
response of the domain through the coupling with elasticity. Panel (b)
shows the horizontal displacement $u_x$. The displacement changes sign
with respect to the vertical centerline of the domain, corresponding to
horizontal motion in opposite directions on the two sides of the region
where the macrophage concentration is largest. Panel (c) shows the
vertical displacement $u_y$, which is larger near the upper part of the
domain and produces an upward deformation of the free boundary.

Panels (d) and (e) show the concentrations of macrophages $M$ and
cytokines $A$, respectively, at $t=0.05$. The macrophage distribution
is modified by the imposed influx through the upper boundary, leading
to larger values of $M$ near this boundary, while the cytokine
concentration remains concentrated mainly around the central region.
Finally, panel (f) shows the resulting configuration of the moving
domain. The initially flat upper boundary becomes curved, with a larger
upward displacement around the central part of the domain. Thus,
Figure~\ref{fig:RSAMD2D} illustrates the two-way coupling between the
concentration fields and the mechanical response: the spatial
distribution of macrophages induces deformation through the elasticity
equations, while the concentrations evolve on the resulting moving
domain.

These results show how the atherosclerotic plaque evolves over time due to the accumulation of inflammatory mediators and the elastic response of the arterial wall.

\subsection{Wave propagation}

In this section, we aim to illustrate the propagation of the solution, which
represents the spread of inflammation in the form of a traveling wave, despite
the evolution of the domain.

We recall that, in \cite{elkhatib2012}, the authors
considered a reaction--diffusion system is a two-dimensional geometry without coupling to elasticity. More
precisely, their model corresponds to the third and fourth equations of
system~\eqref{eq:reducedsystemathero}, together with the associated boundary
conditions. This model builds on the earlier work presented in
\cite{ElKhatib2007}. The main difference between the two formulations lies in
the spatial dimension. In the one-dimensional model in
\cite{ElKhatib2007}, the recruitment of
monocytes is represented by a reaction term, whereas in the two-dimensional
model, this recruitment occurs through the boundary and is therefore introduced
as a nonlinear boundary condition. When the thickness of the two-dimensional
domain is sufficiently small, the system exhibits a behavior similar to that
of the one-dimensional model. In \cite{elkhatib2012}, the authors also studied
traveling-wave solutions connecting the equilibria identified in
\cite{ElKhatib2007}, thus describing the propagation of inflammation along the
arterial wall. According to the
values of their components, these equilibria represent the inflammatory and
non-inflammatory states. The traveling wave connecting these states thus
describes the propagation of inflammation along the arterial wall.

In this section, we aim to study the wave propagation for system~\eqref{eq:reducedsystemathero} under small deformations of the domain. To this
end, we consider a given set of parameters and compute the corresponding
equilibria of the reaction--diffusion system introduced in
\cite{ElKhatib2007}. As an initial condition, we consider a small rectangular
region located near the left edge of the domain, where the concentrations of
$M$ and $A$ are set to the values of the stable inflammatory equilibrium, while
outside this region they are set to the values of the other equilibrium. This
initial configuration is used to initiate the propagation of the inflammatory
state through the domain. The parameters $\gamma_2$ and $\beta_2$, which control
the coupling between the elasticity and reaction--diffusion equations, are chosen
to be sufficiently small. This choice allows us to observe the wave propagation while the domain undergoes a small deformation. For larger values of the coupling parameters, the domain evolves more rapidly, making the effect of its deformation more significant than the propagation of the solution.

\begin{figure} \begin{subfigure}[b]{0.48\textwidth} \centering \includegraphics[width=\textwidth]{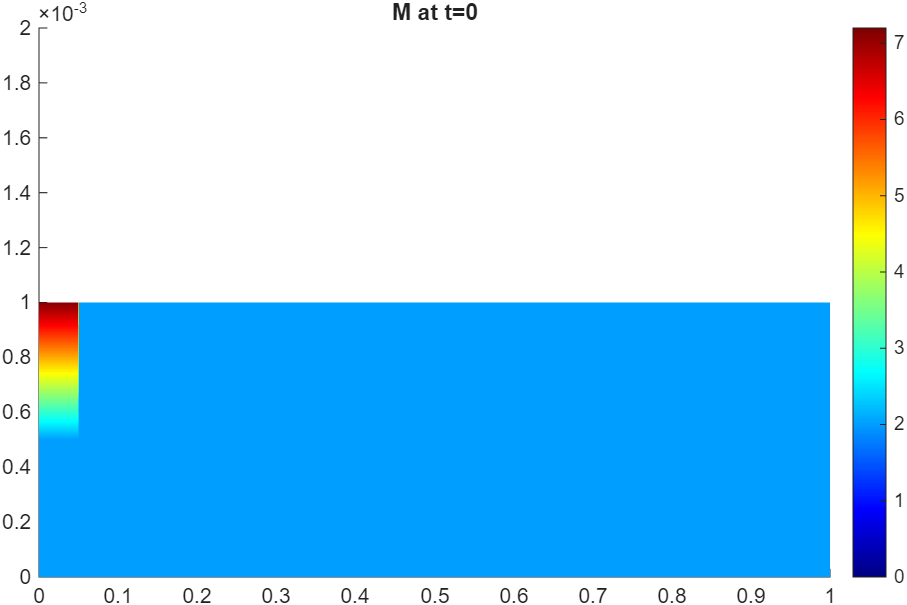} \caption{The concentration of $M$ at time $t=0$} \label{fig:T16_t0} \end{subfigure} \hfill \begin{subfigure}[b]{0.48\textwidth} \centering \includegraphics[width=\textwidth]{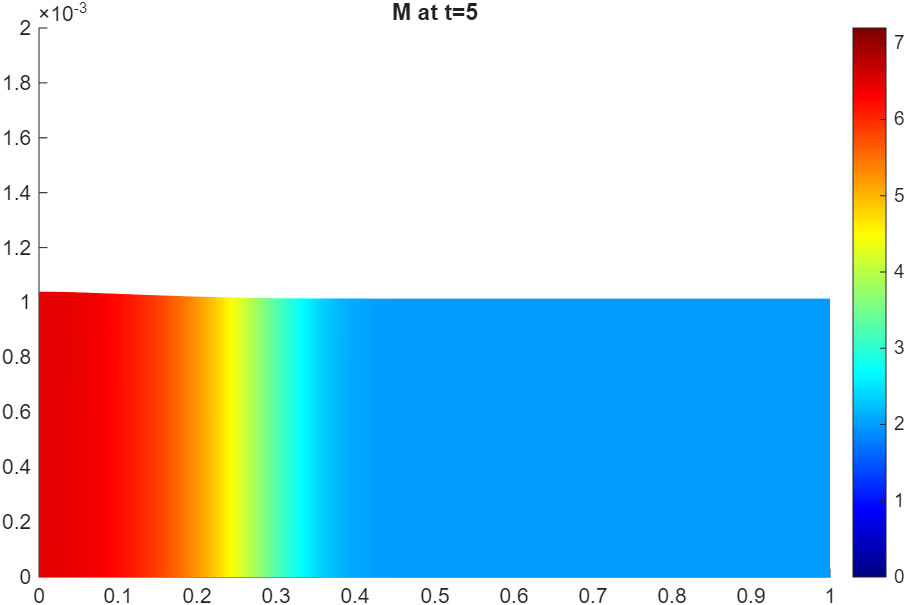} \caption{The concentration of $M$ at time $t=5$} \label{fig:T16_t5} \end{subfigure} \vspace{0.3cm} \begin{subfigure}[b]{0.48\textwidth} \centering \includegraphics[width=\textwidth]{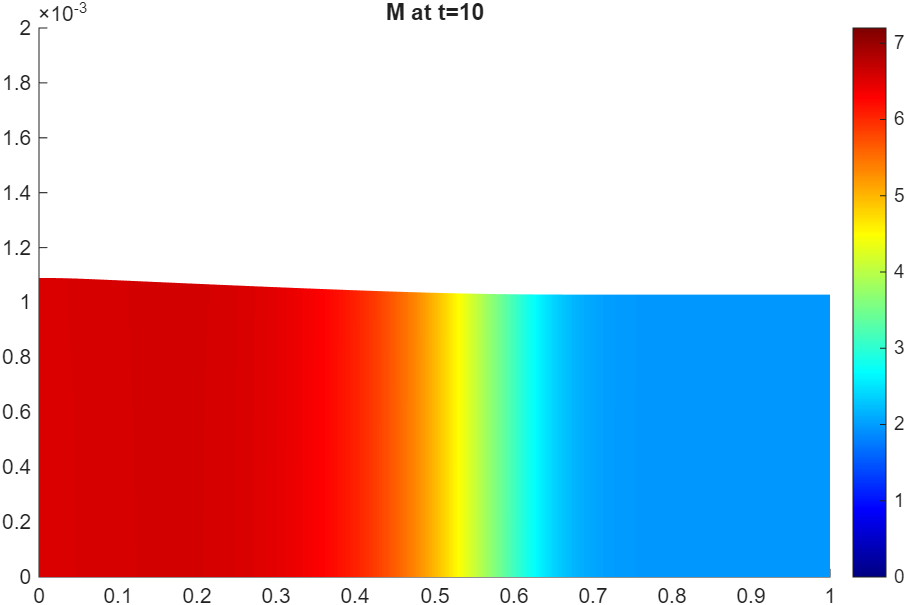} \caption{The concentration of $M$ at time $t=10$} \label{fig:T16_t10} \end{subfigure} \hfill \begin{subfigure}[b]{0.48\textwidth} \centering \includegraphics[width=\textwidth]{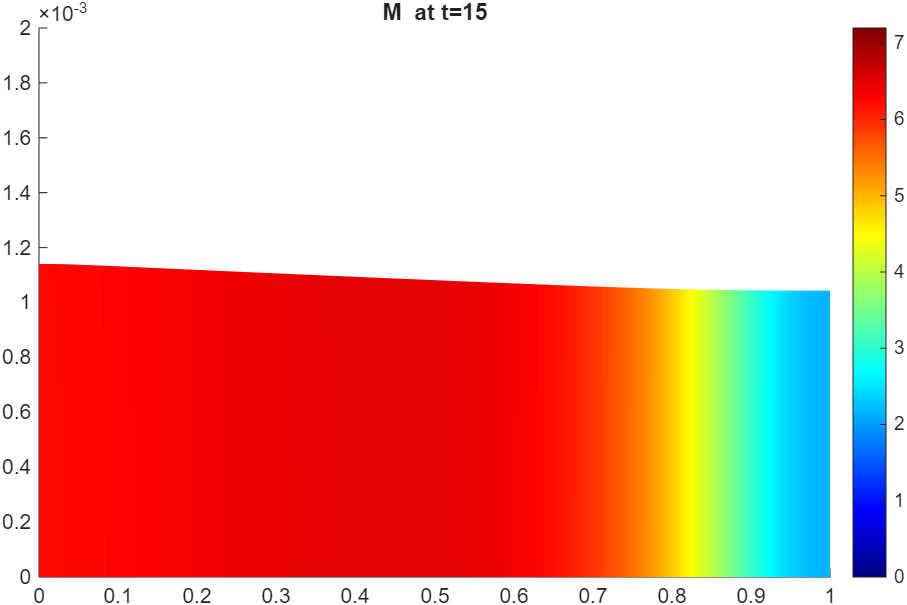} \caption{The concentration of $M$ at time $t=15$} \label{fig:T16_t15} \end{subfigure} \caption{Evolution of the concentration of  $M$ at different times. The parameters taken are
$d_1 = 0.001, \quad d_2 = 0.001, \quad
\alpha_1 = 2, \quad \alpha_2 = 1,\quad
\beta_1 = 8, \quad \lambda_1 = 1, \quad
\lambda_2 = 1, \quad \tau_1 = 1,\quad \tau_2 = 42/43,\quad
\mu = 10, \quad \lambda = 15,\quad\gamma_2 = 10^{-4}, \quad \text{and}\quad \beta_2 = 10^{-4}.
$ The time step is $\Delta t = 0.002$ and the mesh consists of $6004$ elements. The initial concentration of $A$ is taken to be the same as that of $M$.
 } \label{fig:T16_M_evolution} \end{figure}

Figures~\ref{fig:T16_t0}--\ref{fig:T16_t15} show the evolution of the solution
in a rectangular domain of width $1$ and height $10^{-3}$, corresponding to a
small domain thickness. For the set of parameters considered in these simulations,
the model in \cite{ElKhatib2007} admits two equilibria, with $(M,A)=(7.14285,6)$
and $(M,A)=(2,0)$, corresponding to the inflammatory and non-inflammatory states,
respectively.
These values are used as reference states for the two-dimensional
simulations. In the two-dimensional setting, however, the concentrations
may also depend on the transverse coordinate $y$, in particular because
the recruitment of macrophages is prescribed through the upper boundary.
Consequently, the states observed numerically are not necessarily uniform
across the thickness of the domain.
The figures illustrate the propagation of the solution between these
two states, while simultaneously showing the evolution of the upper boundary due
to the coupling with elasticity. We observe that, unlike the results obtained for
the fixed-domain models in \cite{ElKhatib2007} and \cite{elkhatib2012}, the
concentration of $M$ behind the propagating front is slightly lower than its
equilibrium value $7.14285$.

 \begin{figure}[htbp]
    \centering
    \includegraphics[width=0.7\textwidth]{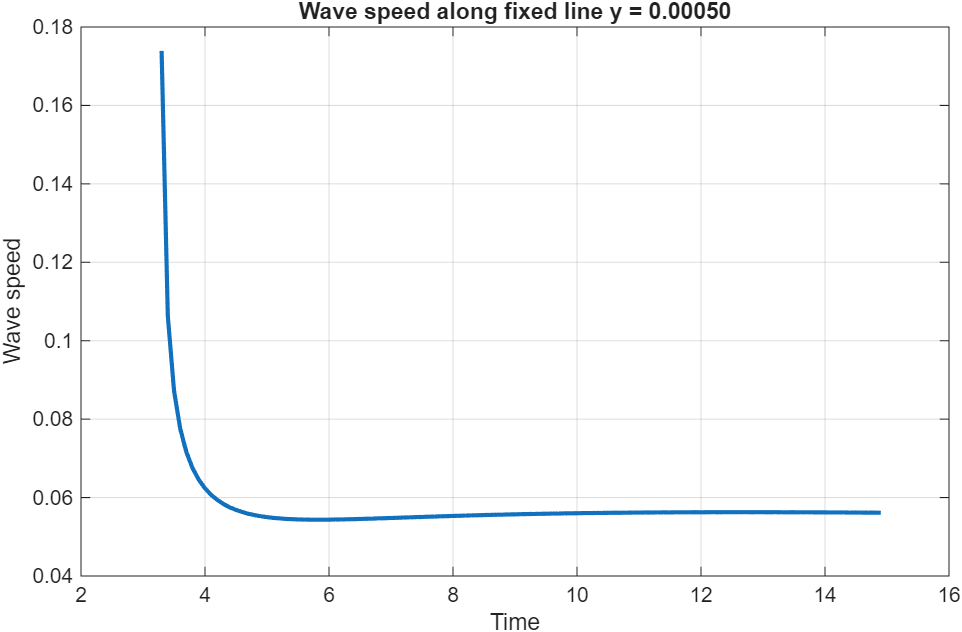}
    \caption{Wave speed of the traveling wave solution computed along the
    horizontal line $y=0.0005$. At each time step, the macrophage concentration $M(x,y,t)$ is evaluated
along the fixed horizontal line $y=0.0005$, and the position of the
propagating front is determined as the value of $x$ at which
$M=4.571425$, corresponding to the midpoint between the inflammatory and
non-inflammatory equilibrium values of $M$. The wave speed is then
computed from the variation of this front position with respect to time.}
    \label{fig:WavespeedT16}
\end{figure}

Figure~\ref{fig:WavespeedT16} shows the evolution of the wave speed computed
along the horizontal line $y=0.0005$. After an initial transient phase, the
wave speed stabilizes and becomes approximately constant, with a value close
to $0.056$. This indicates that the solution exhibits traveling-wave behavior
despite the deformation of the domain.

\subsubsection{Effect of the parameters $\alpha_1$ and $\beta_1$ on the wave speed and domain deformation}

In this subsection, we investigate the effect of the parameters $\alpha_1$ and
$\beta_1$ on the speed of the traveling wave. These parameters are chosen since
they are directly related to the recruitment of monocytes. We recall that
$\alpha_1$ represents the initial recruitment of monocytes associated with the
presence of ox-LDL, while $\beta_1$ characterizes the amplification of monocyte
recruitment induced by inflammatory cytokines.

In the following simulations, $\beta_1$ is varied proportionally to $\alpha_1$.
Indeed, when $\alpha_1=0$, no monocyte influx is considered, and consequently
the contribution associated with $\beta_1$ should also vanish.
Similarly, when $\alpha_1$ increases or decreases, the contribution to the
monocyte influx is assumed to increase or decrease accordingly. Therefore, we
keep a fixed ratio between $\beta_1$ and $\alpha_1$ when studying their effect
on the wave speed.

\begin{figure}[htbp]
    \centering

    \begin{subfigure}[b]{0.48\textwidth}
        \centering
        \includegraphics[width=\textwidth]{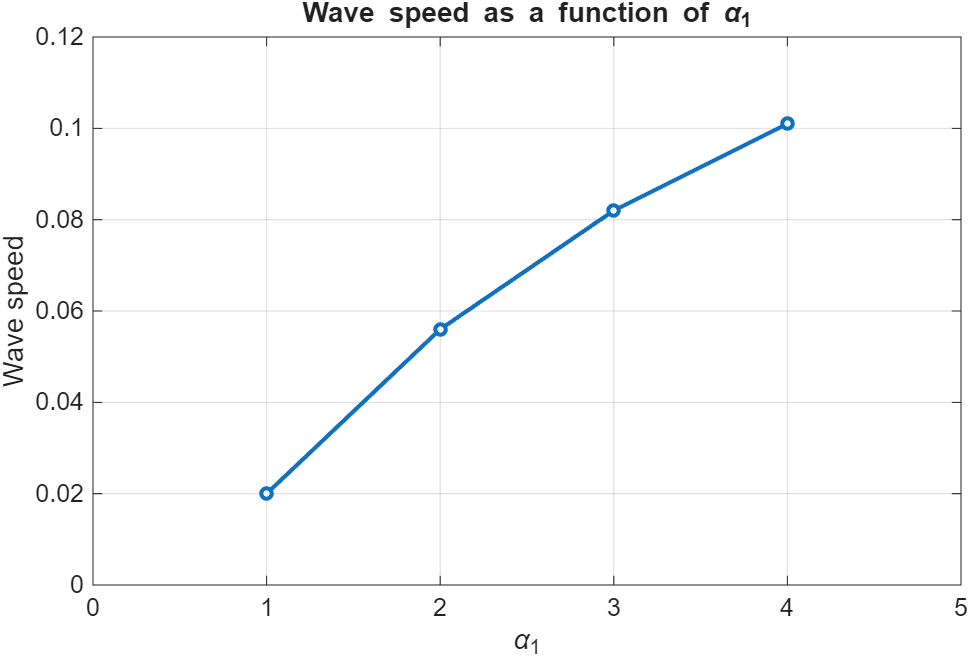}
        \caption{Wave speed of the traveling wave solution as a function of the parameter $\alpha_1$.}
        \label{fig:WavespeedAlpha1}
    \end{subfigure}
    \hfill
    \begin{subfigure}[b]{0.48\textwidth}
        \centering
        \includegraphics[width=\textwidth]{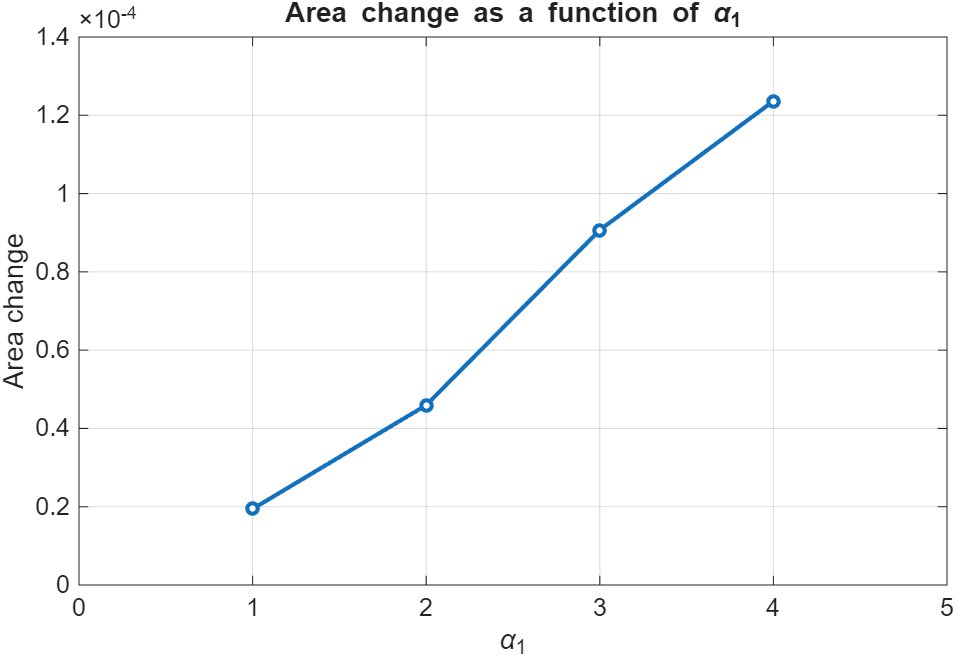}
        \caption{Area change of the domain at time $t=10$ as a function of the parameter $\alpha_1$.}
        \label{fig:AreaChangeAlpha1}
    \end{subfigure}

    \caption{The parameters taken are
    $d_1 = 0.001$, $d_2 = 0.001$, $\alpha_2 = 1$,
    $\beta_1 = 4\alpha_1$, $\lambda_1 = 1$,
    $\lambda_2 = 1$, $\tau_1 = 1$, $\tau_2 = 42/43$,
    $\mu = 10$, $\lambda = 15$, $\gamma_2 = 10^{-4}$,
    and $\beta_2 = 10^{-4}$.}
    \label{fig:WaveSpeedAreaChangeAlpha1}
\end{figure}

In Figure \ref{fig:WaveSpeedAreaChangeAlpha1}, we vary the
parameters $\alpha_1$ and $\beta_1$ while keeping the remaining parameters unchanged and compare the
corresponding wave speeds. Since changing $\alpha_1$ and $\beta_1$ also modifies
the equilibrium states of the reaction--diffusion system, the corresponding
equilibria are computed for each set of parameter values and used to define the
initial conditions. Figure~\ref{fig:WavespeedAlpha1} shows that the wave speed increases as
$\alpha_1$ increases. Since $\alpha_1$ is associated with the influx of
monocytes, this indicates that a higher monocyte influx promotes a faster
propagation of inflammation when the deformation of the domain remains small.
Moreover, Figure~\ref{fig:AreaChangeAlpha1} shows that the change in the domain
area increases with $\alpha_1$, indicating that a higher monocyte influx also
leads to a greater deformation of the domain. These results are consistent with
the biological interpretation, since an increased recruitment of monocytes
enhances the inflammatory response, leading to a faster spread of inflammation
along the arterial wall and a greater deformation of the arterial tissue.

\subsubsection{Effect of the thickness of the domain on the wave speed and domain deformation}
In this section, we investigate the influence of the thickness of the arterial wall on both the propagation speed of the inflammatory wave and the resulting mechanical deformation of the domain. Since the reaction--diffusion processes occur within a deformable medium, changing the thickness of the domain modifies the interaction between the inflammatory dynamics and the elastic response of the arterial wall. To quantify this effect, we compare the wave speed and the area change of the domain for different values of the thickness while keeping all the remaining model parameters fixed.

\begin{figure}[htbp]
    \centering

    \begin{subfigure}[b]{0.48\textwidth}
        \centering
        \includegraphics[width=\textwidth]{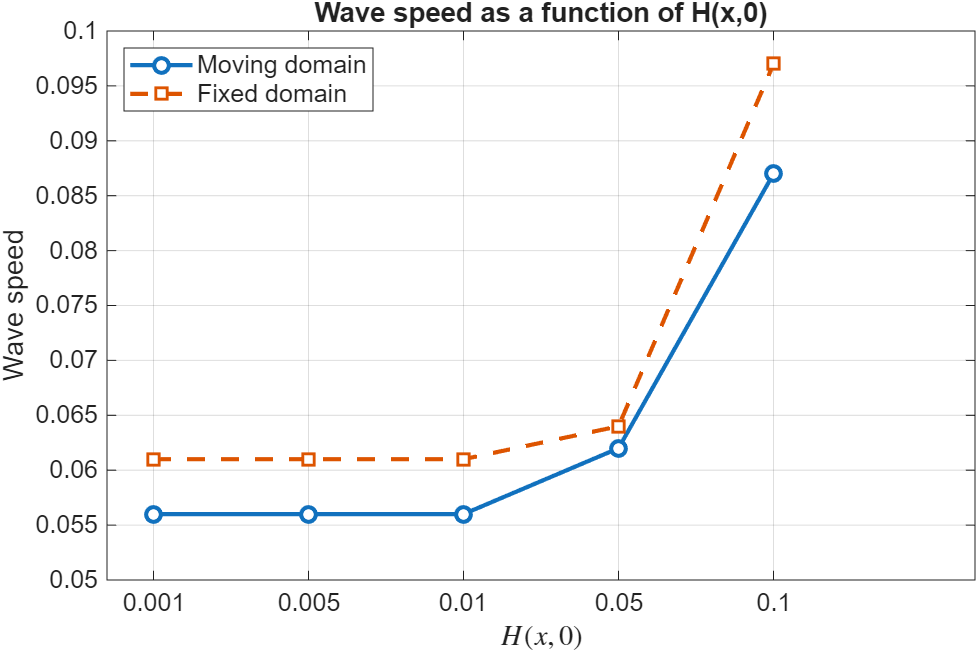}
        \caption{Wave speed of the traveling wave solution as a function of the initial domain thickness $H(x,0)$.}
        \label{fig:WavespeedvsH}
    \end{subfigure}
    \hfill
    \begin{subfigure}[b]{0.48\textwidth}
        \centering
        \includegraphics[width=\textwidth]{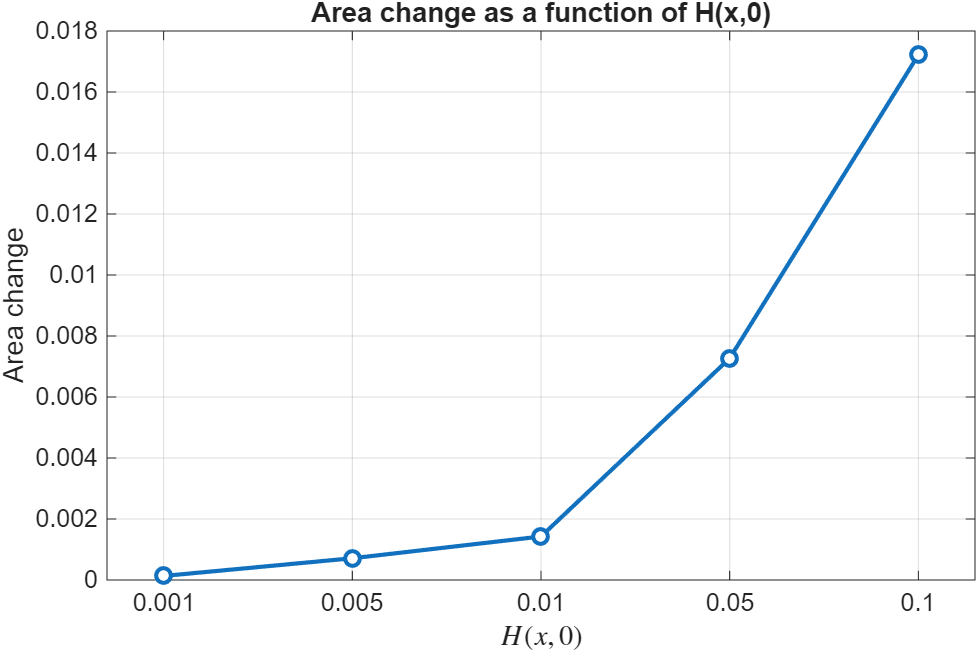}
        \caption{Area change of the domain at time $t=20$ as a function of the initial domain thickness $H(x,0)$.}
        \label{fig:AreaChangeThickness}
    \end{subfigure}

    \caption{The parameters taken are
    $d_1 = 0.001$, $d_2 = 0.001$,
    $\alpha_1 = 2$, $\alpha_2 = 1$,
    $\beta_1 = 8$, $\lambda_1 = 1$,
    $\lambda_2 = 1$, $\tau_1 = 1$, $\tau_2 = 42/43$,
    $\mu = 10$, $\lambda = 15$, $\gamma_2 = 10^{-4}$,
    and $\beta_2 = 10^{-4}$, while the initial domain thickness $H(x,0)$ is varied.}
    \label{fig:WaveSpeedAreaChangeThickness}
\end{figure}

The numerical results are presented in Figure~\ref{fig:WaveSpeedAreaChangeThickness}. In particular, Figure~\ref{fig:WavespeedvsH} compares the wave speed obtained from system~\eqref{eq:reducedsystemathero}, together with the boundary conditions~\eqref{RSABC1}--\eqref{RSABC3}, with the wave speed obtained from the corresponding reaction--diffusion system consisting only of the equations for \(M\) and \(A\) from system~\eqref{eq:reducedsystemathero}, posed on a fixed domain with its corresponding boundary conditions, namely the model studied in~\cite{elkhatib2012}. This comparison is performed for different values of the arterial wall thickness in order to assess the influence of the domain displacement on the propagation speed of the inflammatory wave and on the resulting change in the domain area.

Figure~\ref{fig:WavespeedvsH} shows that, for all values of the thickness considered, the wave speed obtained from system~\eqref{eq:reducedsystemathero} is smaller than the corresponding wave speed of the fixed-domain model of~\cite{elkhatib2012}. Furthermore, both models exhibit the same qualitative behavior as the thickness varies. For relatively small values of the thickness, the wave speed remains nearly constant in both models. As the thickness increases, the wave speed increases in both cases, although it remains consistently lower for the moving-domain model.

These results indicate that increasing the thickness of the arterial wall promotes a faster propagation of the inflammatory wave. However, accounting for the displacement of the domain slows down the propagation compared with the corresponding fixed-domain model. This suggests that the deformation of the arterial wall influences the spread of inflammation, thereby reducing the propagation speed. Consequently, neglecting the domain displacement may lead to an overestimation of the rate at which inflammation propagates through the arterial wall.

Figure~\ref{fig:AreaChangeThickness} shows the variation of the domain area at the fixed time \(t=20\) for different values of the arterial wall thickness. It can be seen that the area change increases as the thickness of the domain increases. This indicates that thicker arterial walls undergo a larger overall deformation under the action of the inflammatory process. From a biological perspective, this suggests that inflammation induces a more pronounced expansion of thicker arterial walls. Combined with the results of Figure~\ref{fig:WavespeedvsH}, these observations show that increasing the thickness of the arterial wall enhances both the propagation speed of the inflammatory wave and the magnitude of the tissue deformation. Moreover, although accounting for the domain displacement reduces the wave speed compared with the corresponding fixed-domain model, the resulting deformation becomes more significant as the wall thickness increases, highlighting the important role of tissue mechanics in the progression of inflammation.

\subsubsection{Effect of the coupling parameter $\gamma_2$
 on the wave speed and domain deformation}
\begin{figure}[htbp]
    \centering

    \begin{subfigure}[b]{0.48\textwidth}
        \centering
        \includegraphics[width=\textwidth]{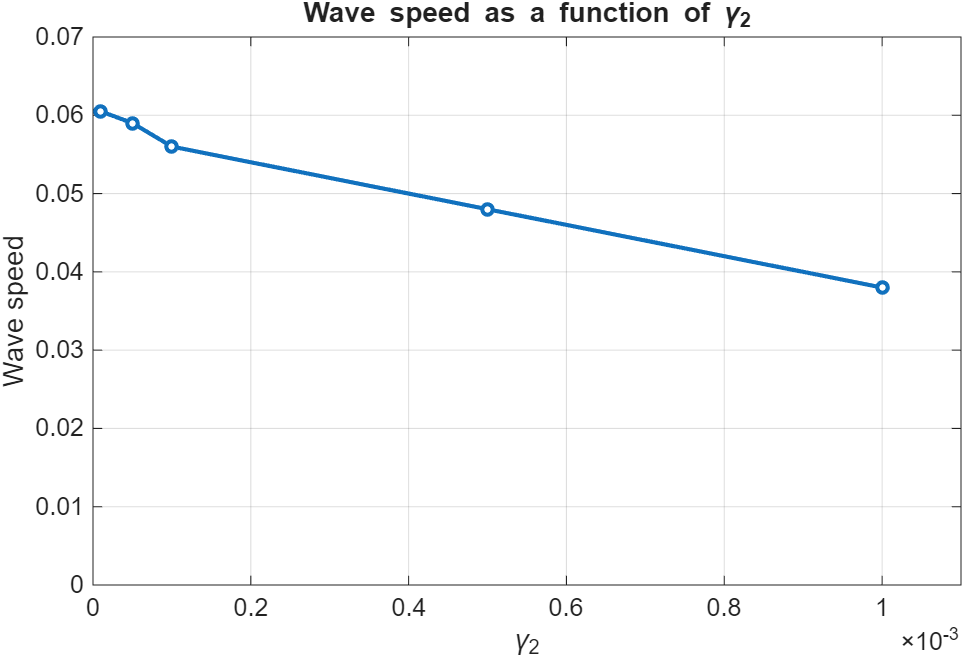}
        \caption{Wave speed of the traveling wave solution as a function of the parameter $\gamma_2$.}
        \label{fig:WavespeedGamma2}
    \end{subfigure}
    \hfill
    \begin{subfigure}[b]{0.48\textwidth}
        \centering
        \includegraphics[width=\textwidth]{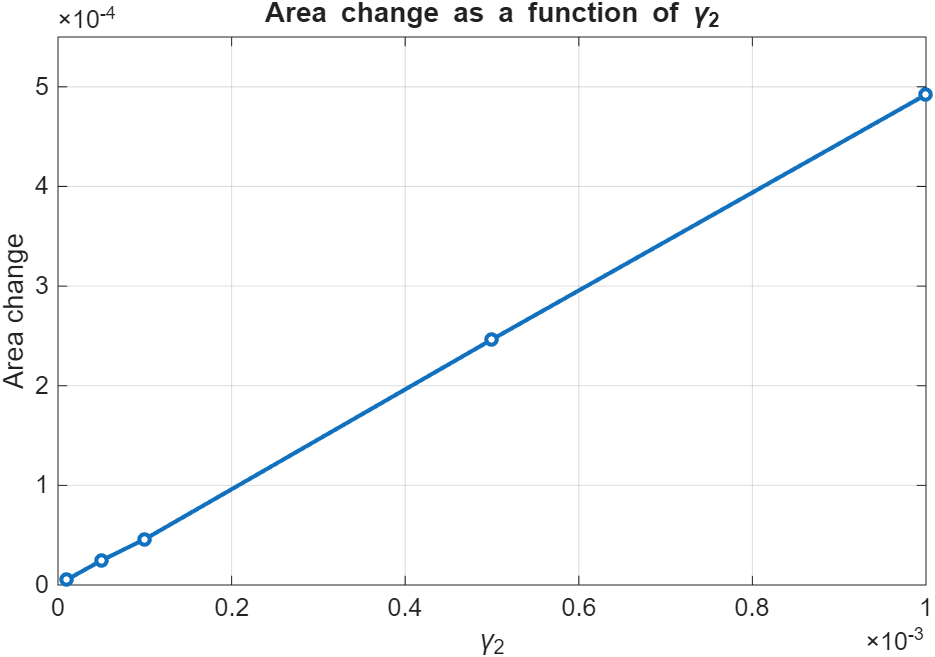}
        \caption{Area change of the domain at time $t=10$ as a function of the parameter $\gamma_2$.}
        \label{fig:AreaChangeGamma2}
    \end{subfigure}

    \caption{The parameters taken are
    $d_1 = 0.001$, $d_2 = 0.001$, $\alpha_1 = 2$,
    $\alpha_2 = 1$, $\beta_1 = 8$, $\lambda_1 = 1$,
    $\lambda_2 = 1$, $\tau_1 = 1$, $\tau_2 = 42/43$,
    $\mu = 10$, $\lambda = 15$,
    and $\beta_2 = \gamma_2$.}
    \label{fig:WaveSpeedAreaChangeGamma2}
\end{figure}

In this section, we investigate the effect of the coupling parameter $\gamma_2$ on both the wave speed and the deformation of the domain. The parameter $\gamma_2$ represents the strength of the coupling between the reaction--diffusion equations. Figure~\ref{fig:WavespeedGamma2} shows that increasing $\gamma_2$ leads to a decrease in the propagation speed of the traveling wave. On the other hand, Figure~\ref{fig:AreaChangeGamma2} shows that the area change induced by the motion of the upper boundary increases as $\gamma_2$ increases. These results indicate that a stronger coupling transfers a larger part of the effect of the inflammatory concentrations to the mechanical deformation of the tissue.
Biologically, this suggests that increasing the coupling parameter $\gamma_2$ causes the accumulation of macrophages to induce a stronger local deformation or thickening of the arterial wall. At the same time, the resulting deformation and enlargement of the domain may redistribute the inflammatory concentrations over a larger region, thereby reducing the spatial propagation speed of the inflammatory solution.
Thus, increasing $\gamma_2$ enhances the mechanical response of the arterial wall while reducing the speed at which inflammation propagates through the tissue.

\subsection{One-dimensional reduction at a fixed horizontal position}

We consider the two-dimensional domain
\[
D(t)
=
\left\{
(x,y)\in\mathbb{R}^2
\;\middle|\;
0\leq x\leq L,\quad
0\leq y\leq H(x,t)
\right\}.
\]
Let \(x=x_0\) be fixed and define
\[
h(t)=H(x_0,t).
\]

We assume that the solution
is independent of \(x\). We therefore write
\[
M=M(y,t),\qquad
A=A(y,t),\qquad
u_x=u_x(y,t),\qquad
u_y=u_y(y,t),
\]
for
\[
0<y<h(t).
\]

Under this assumption, the two-dimensional system reduces to
\begin{equation}
\left\{
\begin{aligned}
\mu \frac{\partial^2u_x}{\partial y^2}
&=0,
\\[4pt]
(\lambda+2\mu)
\frac{\partial^2u_y}{\partial y^2}
-
\gamma_2\frac{\partial M}{\partial y}
&=0,
\\[4pt]
\frac{\partial M}{\partial t}
&=
d_1\frac{\partial^2M}{\partial y^2}
-\lambda_1M,
\\[4pt]
\frac{\partial A}{\partial t}
&=
d_2\frac{\partial^2A}{\partial y^2}
+
\frac{\alpha_2A}{1+A/\tau_2}M
-\lambda_2A.
\end{aligned}
\right.
\label{eq:1D_reduced_system}
\end{equation}

The components of the stress tensor become
\begin{equation}
\begin{aligned}
\sigma_{xx}
&=
\lambda\frac{\partial u_y}{\partial y}
-\beta_2M,
\\[4pt]
\sigma_{xy}
=
\sigma_{yx}
&=
\mu\frac{\partial u_x}{\partial y},
\\[4pt]
\sigma_{yy}
&=
(\lambda+2\mu)\frac{\partial u_y}{\partial y}
-\beta_2M.
\end{aligned}
\label{eq:1D_stress_components}
\end{equation}

At the lower boundary \(y=0\), the boundary conditions are
\begin{equation}
\frac{\partial M}{\partial y}(0,t)=0,
\qquad
\frac{\partial A}{\partial y}(0,t)=0,
\qquad
u_x(0,t)=0,
\qquad
u_y(0,t)=0.
\label{eq:1D_lower_BC}
\end{equation}

In the one-dimensional setting, the outward unit normal at the upper
boundary \(y=h(t)\) is
\[
\boldsymbol n=(0,1).
\]
Consequently, the boundary conditions reduce to
\[
\sigma_{yx}(h(t),t)=0,
\qquad
\sigma_{yy}(h(t),t)=0.
\]
Therefore,
\begin{equation}
\frac{\partial u_x}{\partial y}(h(t),t)=0,
\label{eq:ux_upper_BC}
\end{equation}
and
\begin{equation}
(\lambda+2\mu)
\frac{\partial u_y}{\partial y}(h(t),t)
-
\beta_2M(h(t),t)
=
0.
\label{eq:uy_upper_BC}
\end{equation}

The concentration boundary conditions at \(y=h(t)\) become
\begin{equation}
\frac{\partial M}{\partial y}(h(t),t)
=
\frac{1}{d_1}
\frac{
\alpha_1+\beta_1A(h(t),t)
}{
1+A(h(t),t)/\tau_1
},
\label{eq:M_upper_BC_1D}
\end{equation}
and
\begin{equation}
\frac{\partial A}{\partial y}(h(t),t)=0.
\label{eq:A_upper_BC_1D}
\end{equation}

The complete one-dimensional problem is thus
\begin{equation}
\left\{
\begin{aligned}
\mu u_{x,yy}&=0,
&&0<y<h(t),
\\[3pt]
(\lambda+2\mu)u_{y,yy}-\gamma_2M_y&=0,
&&0<y<h(t),
\\[3pt]
M_t&=d_1M_{yy}-\lambda_1M,
&&0<y<h(t),
\\[3pt]
A_t&=
d_2A_{yy}
+
\frac{\alpha_2A}{1+A/\tau_2}M
-\lambda_2A,
&&0<y<h(t),
\\[5pt]
M_y(0,t)&=0,
\qquad
A_y(0,t)=0,
\\[3pt]
u_x(0,t)&=0,
\qquad
u_y(0,t)=0,
\\[3pt]
u_{x,y}(h(t),t)&=0,
\\[3pt]
(\lambda+2\mu)u_{y,y}(h(t),t)
-\beta_2M(h(t),t)&=0,
\\[3pt]
M_y(h(t),t)
&=
\frac{
\alpha_1+\beta_1A(h(t),t)
}{
d_1\left(1+A(h(t),t)/\tau_1\right)
},
\\[8pt]
A_y(h(t),t)&=0.
\end{aligned}
\right.
\label{eq:complete_1D_problem}
\end{equation}

The initial conditions are
\begin{equation}
M(y,0)=M_0(y),
\qquad
A(y,0)=A_0(y),
\qquad
0\leq y\leq h(0).
\label{eq:1D_initial_conditions}
\end{equation}

The mechanical part of the one-dimensional system can be solved explicitly.
Proceeding similarly to the one-dimensional mechanical problem considered
in Section~4.2, and using the corresponding boundary conditions, we obtain
\[
u_x(y,t)=0.
\]
For the vertical displacement, we obtain
\[
u_y(y,t)
=
\frac{\gamma_2}{\lambda+2\mu}
\int_0^y M(s,t)\,ds
+
\frac{\beta_2-\gamma_2}{\lambda+2\mu}
yM(h(t),t).
\]

In the particular case $\beta_2=\gamma_2$, the displacement of the upper boundary is given by

\[
u_y(h(t),t)
=
\frac{\beta_2}{\lambda+2\mu}
\int_0^{h(t)}M(s,t)\,ds.
\]

This last relation provides a direct link between the mechanical deformation
of the arterial wall and the total quantity of macrophages within the domain.
Indeed, the integral
\[
\int_0^{h(t)} M(s,t)\,ds
\]
represents the total macrophage concentration across the thickness of the
wall. Therefore, when $\beta_2=\gamma_2$, the displacement of the upper
boundary, and hence the change in wall thickness, is proportional to the
total quantity of macrophages accumulated in the domain. Thus, an increase in macrophage
accumulation results in a corresponding increase in the wall thickness.

\section{Discussion}

The present work develops a coupled elasticity--diffusion framework for describing the interaction between diffusion and mechanical deformation in deformable media. Two numerical approaches are proposed for the two-dimensional formulation, corresponding to fixed and moving domain computations, allowing the influence of domain evolution on the solution to be investigated. A quasi-stationary formulation is then derived together with corresponding numerical simulations, providing a simplified framework for studying the coupled system. In addition, a linear mass--volume relation is established, showing that, under the assumptions of linear elasticity and small deformations, the total mass of the diffusing species is proportional to the volume change of the domain. A one-dimensional elasticity--diffusion model is also introduced as a particular case of the two-dimensional formulation, where analytical properties of the coupled system are investigated. Finally, the framework is applied to a model of atherosclerosis. Numerical simulations are used to study inflammatory wave propagation and to quantify the influence of the mechanical coupling and arterial wall thickness on the wave speed and arterial wall deformation. A one-dimensional reduction at a fixed horizontal position is further derived, providing a simplified description of the interaction between inflammation and the local deformation of the arterial wall.

Overall, the analytical and numerical results show that coupling diffusion with elasticity creates a reciprocal interaction between the concentration fields and the geometry of the domain. A nonuniform concentration generates horizontal and vertical displacements, with the largest deformation occurring in regions of higher concentration, while increasing the mechanical coupling strengthens the deformation. The comparison of the fixed and moving domain approaches indicates that both formulations give similar predictions for small displacements, whereas their results progressively differ as the deformation increases because the moving domain formulation accounts for the evolving geometry. The derived mass--volume relation further shows that, within the linear-elasticity and small-deformation regime, the change in domain volume is proportional to the total mass of the diffusing species. The numerical results confirm this relation when the displacement remains small, while noticeable deviations appear for larger deformations, for which the neglected higher-order geometric terms are no longer negligible. The one-dimensional formulation and the reduction of the two-dimensional problem at a fixed horizontal position also provide an explicit relation between the displacement of the moving boundary and the total concentration across the thickness of the domain, thereby linking local species accumulation directly to wall thickening.

In the application to atherosclerosis, the model reproduces the propagation of an inflammatory front despite the simultaneous deformation of the arterial wall, and the wave speed approaches an approximately constant value after an initial transient. The macrophage concentration behind the front is slightly lower than the corresponding equilibrium value of the fixed-domain model, indicating that the evolving geometry modifies the inflammatory state. The parameter studies show that stronger monocyte recruitment accelerates the propagation of inflammation and increases wall deformation, while increasing the arterial-wall thickness enhances both the wave speed and the overall domain expansion. Nevertheless, the inflammatory wave propagates more slowly in the moving-domain formulation than in the corresponding fixed-domain model, suggesting that neglecting wall deformation may overestimate the rate of inflammation propagation. Moreover, increasing the strength of the mechanical coupling increases the domain deformation while reducing the wave speed, which indicates that a larger part of the inflammatory effect is transferred to tissue expansion and that the resulting enlargement of the domain may redistribute the concentrations and slow their spatial propagation.

These results should be interpreted in light of the distinction between the present formulation and morphoelastic models of tissue growth. In morphoelastic formulations, tissue growth modifies the local stress-free reference configuration, generally through a multiplicative decomposition of the deformation gradient into elastic and growth components, with the growth tensor prescribed by a constitutive law depending on biochemical concentrations and, possibly, mechanical stresses. In the present model, by contrast, no separate growth tensor or remodelling law is introduced. The macrophage concentration is coupled directly to the quasi-static elasticity equations, and the resulting displacement determines the evolution of the domain boundary. Accordingly, the wall thickening observed in our simulations represents concentration-induced elastic deformation rather than permanent tissue growth or remodelling. Extending the present formulation to include a morphoelastic decomposition and an appropriate growth law constitutes a natural direction for future work.

The relation obtained between macrophage accumulation and boundary displacement also clarifies an important distinction between the present formulation and classical Stefan-type free-boundary models. In a Stefan problem, the normal velocity of the interface is typically determined by the flux of a diffusing quantity across that interface. Here, however, macrophages may migrate and redistribute within the intima, and their accumulation generates mechanical deformation throughout the tissue. Consequently, the arterial-wall boundary can be displaced by cells located inside the domain even when the corresponding mechanism cannot be represented solely by a cellular flux across the moving interface. Our results therefore support a mechanically driven description of boundary evolution, in which cellular accumulation acts through the elastic displacement of the arterial wall.

The numerical treatment of the moving domain may also be viewed in relation to established approaches for evolving-domain problems, including ALE methods ~\cite{donea1982,donea2004}, level-set methods~\cite{osher1988,sethian2003}, and phase-field formulations~\cite{karma1998}. In the present work, the computational mesh is updated directly from the displacement field, resulting in a Lagrangian treatment of domain evolution. The comparison with the fixed-domain formulation shows that the simpler approximation remains accurate when the displacement is small, but progressively loses accuracy as the deformation increases. This observation provides a practical indication of the regime in which solving the equations on a fixed reference configuration may be sufficient and of when the geometric evolution of the domain must be explicitly taken into account.

The present model has also several limitations. First, the mechanical behavior of the arterial wall is described by linear elasticity under the assumption of small deformations. While this approximation is appropriate for the present study, arterial tissues are known to exhibit nonlinear, anisotropic, and heterogeneous mechanical properties that may become significant under larger deformations or at more advanced stages of the disease. Second, the biological model remains intentionally simplified and focuses on the interaction between a limited number of variables. Important mechanisms involved in the development of atherosclerosis, including LDL transport and oxidation, foam-cell formation, smooth muscle cell migration, extracellular matrix remodeling, calcification, and plaque rupture, are not incorporated in the current formulation. Finally, although the model accounts for the deformation of the arterial wall, it does not include the effect of blood flow and therefore neglects the fluid--structure interaction between the hemodynamics and the vessel wall. Incorporating more realistic constitutive laws, additional biological processes, and a fully coupled fluid--structure interaction model constitutes a natural direction for future work and would provide a more comprehensive description of atherosclerosis progression.
\section*{Acknowledgements}

One of the authors (NB) was supported within the framework of the state assignment of the Institute for Problems in Mechanical Engineering of the RAS № 124041500009-8 with the support of the Ministry of Science and Higher Education of the RF.\\
Two of the authors (NK) and (GA) were supported by the PIRF project No I0028 of the Lebanese
American University.
\newpage
\appendix
\section*{Appendix}
\addcontentsline{toc}{section}{Appendix}
\section{Derivation of the Mass--Volume Relation}
\label{appendix:massvolume}

This appendix provides the detailed derivation of the mass--volume relation presented in Section~\ref{sec:massvolume}. The proof is based on integrating the elasticity equations over the domain and using the stress-free boundary conditions.

We begin with the first elasticity equation \eqref{eq:2.7} in the rectangular domain
\[
0<x<L,\qquad 0<y<H,
\]
where $(u_x,u_y)=(u,v)$ and $\gamma_2=\beta_2=\beta$.

Integrating equation \eqref{eq:2.7} with respect to $x$ gives

\begin{align}
\beta(C(x,y)-C(0,y))
&=
(\lambda+2\mu)(u'_x(x,y)-u'_x(0,y))
\nonumber\\
&\quad
+\mu\int_0^x u''_{y^2}(\xi,y)\,d\xi
+
(\lambda+\mu)(v'_y(x,y)-v'_y(0,y)).
\label{eq:2.11}
\end{align}

We now integrate this equation over the whole domain. The different terms become

\begin{align}
\int_0^H\int_0^L u'_x(x,y)\,dx\,dy
&=
\int_0^H
\left(
u(L,y)-u(0,y)
\right)\,dy,
\label{eq:2.12}
\\
\int_0^H\int_0^L u'_x(0,y)\,dx\,dy
&=
L
\int_0^H
u'_x(0,y)\,dy,
\label{eq:2.13}
\\
\int_0^H\int_0^L
\int_0^x
u''_{y^2}(\xi,y)\,d\xi\,dx\,dy
&=
\int_0^L
\int_0^x
\left(
u'_y(\xi,H)-u'_y(\xi,0)
\right)
\,d\xi\,dx,
\label{eq:2.16}
\\
\int_0^H\int_0^L
v'_y(x,y)\,dx\,dy
&=
\int_0^L
\left(
v(x,H)-v(x,0)
\right)\,dx,
\label{eq:2.14}
\\
\int_0^L\int_0^H
v'_y(0,y)\,dx\,dy
&=
L
\left(
v(0,H)-v(0,0)
\right).
\label{eq:2.15}
\end{align}

The stress-free boundary conditions are

\[
\begin{cases}
x=0:&\sigma_{xx}=0,\\
y=0:&\sigma_{xy}=0,\\
y=H:&\sigma_{xy}=0.
\end{cases}
\]

or equivalently,

\begin{align}
x=0:\qquad
&
2\mu u'_x(0,y)
+
\lambda
\left(
u'_x(0,y)+v'_y(0,y)
\right)
=
\beta C(0,y),
\label{bc1}
\\
y=0:\qquad
&
u'_y(x,0)
=
-
v'_x(x,0),
\label{bc2}
\\
y=H:\qquad
&
u'_y(x,H)
=
-
v'_x(x,H).
\label{bc3}
\end{align}

Using \eqref{bc1}, equation \eqref{eq:2.11} becomes

\begin{align}
\beta C(x,y)
&=
(\lambda+2\mu)u'_x(x,y)
+
\mu
\int_0^x
u''_{y^2}(\xi,y)\,d\xi
\nonumber\\
&\quad
+
(\lambda+\mu)
v'_y(x,y)
-
\mu
v'_y(0,y).
\label{eq:2.19}
\end{align}

Next, using the tangential stress conditions \eqref{bc2}--\eqref{bc3}, the integral in \eqref{eq:2.16} can be rewritten as

\begin{align}
\int_0^L
\int_0^x
(u'_y(\xi,H)-u'_y(\xi,0))
\,d\xi\,dx
&=
-
\int_0^L
\int_0^x
(v'_x(\xi,H)-v'_x(\xi,0))
\,d\xi\,dx
\nonumber\\
&=
-
\int_0^L
\Big[
(v(x,H)-v(0,H))
\nonumber\\
&\qquad\qquad
-
(v(x,0)-v(0,0))
\Big]
dx
\nonumber\\
&=
\int_0^L
(v(x,0)-v(x,H))
dx
\nonumber\\
&\quad
+
L(v(0,H)-v(0,0)).
\label{eq:2.20}
\end{align}

Substituting the previous identities into equation \eqref{eq:2.19} yields

\begin{align}
\beta
\int_0^H
\int_0^L
C(x,y)\,dx\,dy
&=
(2\mu+\lambda)
\int_0^H
(u(L,y)-u(0,y))
\,dy
\nonumber\\
&\quad
+
\mu
\int_0^L
(v(x,0)-v(x,H))
dx
\nonumber\\
&\quad
+
\mu
L
(v(0,H)-v(0,0))
\nonumber\\
&\quad
+
(\lambda+\mu)
\int_0^L
(v(x,H)-v(x,0))
dx
\nonumber\\
&\quad
-
\mu
L
(v(0,H)-v(0,0))
\nonumber\\
&=
(2\mu+\lambda)
\int_0^H
(u(L,y)-u(0,y))
\,dy
\nonumber\\
&\quad
+
\lambda
\int_0^L
(v(x,H)-v(x,0))
\,dx.
\label{eq:mass11}
\end{align}

Equation \eqref{eq:mass11} is the first identity used in the derivation of the mass--volume relation.

We now perform the same calculation for equation (\ref{eq:2.8}). Consider the equation in a rectangular domain $0<x<L$, $0<y<H$, with $u_x=u$, $u_y=v$, $\gamma_2=\beta_2=\beta$, it reads:
\[
\beta \frac{\partial C}{\partial y}
=
\mu\Delta v
+
(\lambda+\mu)
\left(
\frac{\partial^2 v}{\partial y^2}
+
\frac{\partial^2 u}{\partial x\partial y}
\right).
\]

Integrating with respect to $y$, we obtain

\begin{align}
\beta(C(x,y)-C(x,0))
&=
(\lambda+2\mu)
(v'_y(x,y)-v'_y(x,0))
\nonumber\\
&\quad
+
\mu
\int_0^y
v''_{x^2}(x,\xi)\, d\xi
+
(\lambda+\mu)
(u'_x(x,y)-u'_x(x,0)).
\label{eq:bis2.11}
\end{align}

Now we integrate with respect to $x$ and $y$ over the whole domain:

\begin{align}
\int_0^H
\int_0^L
v'_y(x,y)\,dx\,dy
&=
\int_0^L
(v(x,H)-v(x,0))\,dx,
\label{eq:bis2.14}
\\
\int_0^H
\int_0^L
v'_y(x,0)\,dx\,dy
&=
H
\int_0^L
v'_y(x,0)\,dx,
\label{eq:bis2.15}
\\
\int_0^H
\int_0^L
\int_0^y
v''_{x^2}(x,\xi)\, d\xi\,dx\,dy
&=
\int_0^H
\int_0^y
(v'_x(L,\xi)-v'_x(0,\xi))
\,d\xi\,dy,
\label{eq:bis2.16}
\\
\int_0^H
\int_0^L
u'_x(x,y)\,dx\,dy
&=
\int_0^H
(u(L,y)-u(0,y))\,dy,
\label{eq:bis2.12}
\\
\int_0^H
\int_0^L
u'_x(x,0)\,dx\,dy
&=
\int_0^H
(u(L,0)-u(0,0))\,dy
=
H(u(L,0)-u(0,0)).
\label{eq:bis2.13}
\end{align}

We consider the following boundary conditions:
\[
\begin{cases}
y=0: & \sigma_{yy}=0,\\
x=0: & \sigma_{xy}=0,\\
x=L: & \sigma_{xy}=0.
\end{cases}
\]

That is,

\begin{align}
y=0:\qquad
&
2\mu v'_y(x,0)
+
\lambda
(u'_x(x,0)+v'_y(x,0))
=
\beta C(x,0),
\label{bc4}
\\
x=0:\qquad
&
u'_y(0,y)
=
-
v'_x(0,y),
\label{bc5}
\\
x=L:\qquad
&
u'_y(L,y)
=
-
v'_x(L,y).
\label{bc6}
\end{align}

Using equation \eqref{bc4}, equation \eqref{eq:bis2.11} becomes

\begin{align}
\beta C(x,y)
&=
(\lambda+2\mu)
v'_y(x,y)
+
\mu
\int_0^y
v''_{x^2}(x,\xi)\,d\xi
\nonumber\\
&\quad
+
(\lambda+\mu)
u'_x(x,y)
-
\mu
u'_x(x,0).
\label{eq:bis2.19}
\end{align}

Using the boundary conditions at $x=0$ and $x=L$, we write
\eqref{eq:bis2.16} as

\begin{align}
\int_0^H
\int_0^y
(v'_x(L,\xi)-v'_x(0,\xi))
\,d\xi\,dy
&=
-
\int_0^H
\int_0^y
(u'_y(L,\xi)-u'_y(0,\xi))
\,d\xi\,dy
\nonumber\\
&=
-
\int_0^H
\left[
(u(L,y)-u(L,0))
-
(u(0,y)-u(0,0))
\right]
dy
\nonumber\\
&=
\int_0^H
(u(0,y)-u(L,y))
dy
+
H(u(L,0)-u(0,0)).
\label{eq:bis2.20}
\end{align}

Thus, after integration of \eqref{eq:bis2.19},

\begin{align}
\beta
\int_0^L
\int_0^H
C(x,y)\,dx\,dy
&=
(2\mu+\lambda)
\int_0^L
(v(x,H)-v(x,0))
\,dx
\nonumber\\
&\quad
+
\mu
\int_0^H
(u(0,y)-u(L,y))
dy
+
\mu
H
(u(L,0)-u(0,0))
\nonumber\\
&\quad
+
(\lambda+\mu)
\int_0^H
(u(L,y)-u(0,y))
\,dy
-
\mu
H
(u(L,0)-u(0,0))
\nonumber\\
&=
\lambda
\int_0^H
(u(L,y)-u(0,y))
\,dy
+
(\lambda+2\mu)
\int_0^L
(v(x,H)-v(x,0))
\,dx.
\label{eq:bis2.21}
\end{align}

From equations \eqref{eq:mass11} and \eqref{eq:bis2.21}, we have
\begin{align}
\beta \int_0^L \int_0^H C(x,y)\,dx\,dy
&= (2\mu+\lambda)\int_0^H (u(L,y)-u(0,y))\,dy
+ \lambda\int_0^L (v(x,H)-v(x,0))\,dx \notag\\
&= \lambda\int_0^H (u(L,y)-u(0,y))\,dy
+ (\lambda+2\mu)\int_0^L (v(x,H)-v(x,0))\,dx.
\end{align}

Therefore, we deduce

\begin{equation}
\int_0^H (u(L,y)-u(0,y))\,dy
=
\int_0^L (v(x,H)-v(x,0))\,dx.
\label{5.28}
\end{equation}

Since the integrals represent the volume change, this means that the total volume change in the horizontal direction is the same as in the vertical direction.

\subsection*{Volume Change}

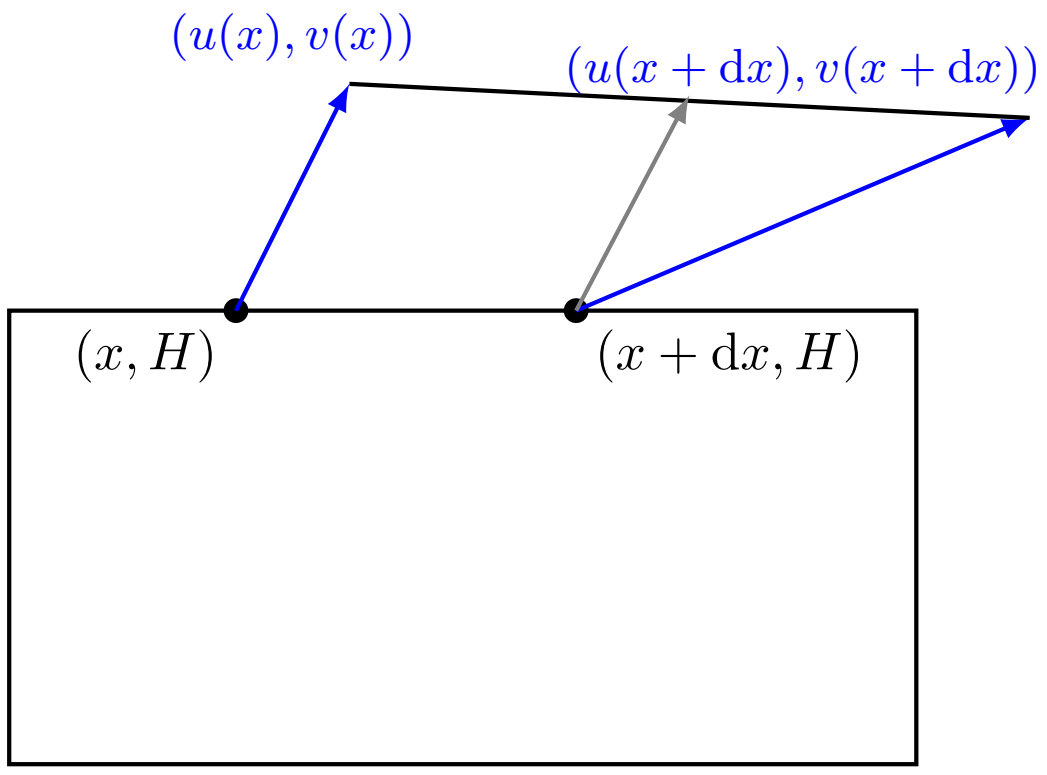
\begin{figure}[h]
\begin{center}
\begin{tikzpicture}[scale=1.5, >=latex]

\def\x{2}
\def\dx{1.5}
\def\H{2}

\def\uone{0.5}
\def\vone{1}

\def\utwo{2}
\def\vtwo{0.85}

\draw[thick] (1,0) rectangle (5,2);

\filldraw[black] (\x,\H) circle (0.05) node[below left] {$(x,H)$};
\filldraw[black] ({\x+\dx},\H) circle (0.05) node[below right] {$(x+\mathrm{d}x,H)$};

\draw[->,thick,blue] (\x,\H) -- ++(\uone,\vone)
node[midway,left=2.5pt,above=22pt]
{\small$(u(x),v(x))$};

\draw[->,thick,blue] ({\x+\dx},\H) -- ++(\utwo,\vtwo)
node[midway,right=40pt,above=19pt]
{\small$(u(x+\mathrm{d}x),v(x+\mathrm{d}x))$};

\draw[thick]
({\x+\uone},{\H+\vone})
--
({\x+\dx+\utwo},{\H+\vtwo});

\draw[->,thick,gray]
({\x+\dx},\H)
--
++(\uone,0.95);

\end{tikzpicture}
\end{center}
\caption{Area of the additional domain.}
\label{volume}
\end{figure}

Let the displacement of $(x,H)$ be $(u(x,H),v(x,H))$ and that of $(x+dx,H)$ be $(u(x+dx,H),v(x+dx,H))$. We need to find the area of the additional domain. It is composed of the parallelogram and triangle. See Figure~\ref{volume}.

The area of the parallelogram is
\[
v(x,H)\,dx,
\]
and the area of the triangle is
\[
\frac{1}{2}(u(x+dx,H)-u(x,H))\,v(x+dx,H)
=
\frac{1}{2}
\frac{u(x+dx,H)-u(x,H)}{dx}
v(x+dx,H)\,dx.
\]

Passing to the continuous limit,

\[
V
=
\int_0^L v(x,H)\,dx
+
\frac12\int_0^L v(x,H)u'(x,H)\,dx.
\]

For small displacements, the second (quadratic) term is negligible, so the total mass in \eqref{eq:mass11} is expressed in terms of the total volume change on all four boundaries of the rectangle.

\subsection*{The Mass--Volume Relation}

The previous relations show that the integral of the concentration can be expressed in terms of boundary displacements. Since these boundary integrals represent the global deformation of the domain, they can be interpreted as the total volume change. From the previous calculations, we conclude that, for small displacements, the total mass of the diffusing species is proportional to the total volume change, denoted by $\Delta V$, of the domain. More precisely, we obtain from equations \eqref{eq:mass11} and \eqref{5.28}

\begin{align}
\Delta V
&=
\int_0^H (u(L,y)-u(0,y))\,dy
+
\int_0^L (v(x,H)-v(x,0))\,dx
\notag\\
&=
2\int_0^L (v(x,H)-v(x,0))\,dx
\notag\\
&=
\frac{\beta}{\lambda+\mu}
\int_0^L \int_0^H C(x,y)\,dx\,dy.
\label{linearc}
\end{align}

This shows that the total volume change is proportional to the total mass, with proportionality constant
\[
\frac{\beta}{\lambda+\mu}.
\]

\section{A geometry-independent form of the mass--volume relation
for the coupled elasticity--diffusion model}
\label{app:general_mass_volume}

\subsection*{Statement of the problem}

In the quasi-stationary formulation of the coupled elasticity--diffusion model,
the mechanical part reduces (with the symmetric coupling
$\gamma_2=\beta_2=\beta$) to the equilibrium of a linearly elastic medium
subject to a concentration-induced pre-stress. Let
$\Om\subset\mathbb{R}^2$ be a bounded domain with Lipschitz boundary $\dOm$
and outward unit normal $\bn$, let $\bu=(u_x,u_y)$ be the displacement field,
$\bveps=\tfrac12(\grad\bu+\grad\bu^{\top})$ the linearized strain, and $C$ the
concentration of the diffusing species. The stress tensor is
\begin{equation}\label{eq:constitutive}
  \bsig \;=\; 2\mu\,\bveps + \lambda\,(\divv\bu)\,\bm{I} - \beta\,C\,\bm{I},
\end{equation}
where $\lambda,\mu$ are the Lam\'e coefficients and $\operatorname{div}\mathbf{u}
=\partial_x u_x+\partial_y u_y$. Written componentwise,
\eqref{eq:constitutive} coincides with the constitutive law
$\sigma_{xx}=2\mu\,\partial_x u_x+\lambda(\partial_x u_x+\partial_y u_y)-\beta C$,
$\sigma_{yy}=2\mu\,\partial_y u_y+\lambda(\partial_x u_x+\partial_y u_y)-\beta C$,
$\sigma_{xy}=\mu(\partial_y u_x+\partial_x u_y)$ of the paper. With this
definition the reduced momentum balance
\[
  \mu\,\Delta u_x+(\lambda+\mu)\,\partial_x(\divv\bu)=\beta\,\partial_x C,
  \qquad
  \mu\,\Delta u_y+(\lambda+\mu)\,\partial_y(\divv\bu)=\beta\,\partial_y C
\]
is exactly the divergence-free condition
\begin{equation}\label{eq:equilibrium}
  \divv\bsig=\bm{0}\quad\text{in }\Om ,
\end{equation}
i.e.\ the concentration coupling appears as a body force that has been absorbed
into $\bsig$ through the pre-stress term $-\beta C\,\bm{I}$.

We assume the medium is entirely traction-free,
\begin{equation}\label{eq:tractionfree}
  \bsig\,\cdot\bn=\bm{0}\qquad\text{on }\dOm .
\end{equation}
For small displacements the change of area (volume, in the plane strain sense)
of the domain is, to first order,
\begin{equation}\label{eq:dV}
  \Delta V=\int_{\Om}\divv\bu\dd\bx
          =\oint_{\dOm}\bu\cdot\bn\dd s .
\end{equation}

\subsection*{Main result}

\begin{theorem}[Mass--volume relation on an arbitrary domain]\label{thm:main}
Let $\Om\subset\mathbb{R}^2$ be a bounded domain with Lipschitz boundary, and
let $(\bu,C)$ satisfy the equilibrium equation \eqref{eq:equilibrium} with the
constitutive law \eqref{eq:constitutive} and the traction-free boundary
condition \eqref{eq:tractionfree}. Then, in the small-deformation regime,
\begin{equation}\label{eq:massvolume}
  \;\Delta V=\frac{\beta}{\lambda+\mu}\int_{\Om}C\dd\bx\;
\end{equation}
independently of the shape of $\Om$.
\end{theorem}

\noindent
The relation \eqref{eq:massvolume} was established in the paper by explicit
integration over a rectangle. Theorem~\ref{thm:main} shows that the rectangular
geometry is inessential: the identity is a direct consequence of the
equilibrium equation and the traction-free condition, and the proof below
requires no assumption on the shape of the domain.

\begin{proof}
Taking the trace of \eqref{eq:constitutive} in two dimensions
($\tr\bm{I}=2$, $\tr\bveps=\divv\bu$) gives
\begin{equation}\label{eq:trace}
  \tr\bsig=\sigma_{xx}+\sigma_{yy}
          =2(\lambda+\mu)\,\divv\bu-2\beta\,C,
\end{equation}
so that
\begin{equation}\label{eq:divu}
  \divv\bu=\frac{\tr\bsig}{2(\lambda+\mu)}+\frac{\beta}{\lambda+\mu}\,C .
\end{equation}
Integrating \eqref{eq:divu} over $\Om$ and using \eqref{eq:dV},
\begin{equation}\label{eq:dVsplit}
  \Delta V=\frac{1}{2(\lambda+\mu)}\int_{\Om}\tr\bsig\dd\bx
          +\frac{\beta}{\lambda+\mu}\int_{\Om}C\dd\bx .
\end{equation}
It remains to show that the first term vanishes, i.e.\ that
$\int_{\Om}\tr\bsig\dd\bx=0$.

Because $\bsig$ is divergence-free, $\partial_k\sigma_{ik}=0$, and therefore for
each pair of indices $i,j$,
\begin{equation}\label{eq:identity}
  \partial_k\bigl(x_j\,\sigma_{ik}\bigr)
   =\delta_{jk}\,\sigma_{ik}+x_j\,\partial_k\sigma_{ik}
   =\sigma_{ij},
\end{equation}
where summation over $k$ is understood. Integrating \eqref{eq:identity} over
$\Om$ and applying the divergence theorem,
\begin{equation}\label{eq:avgstress}
  \int_{\Om}\sigma_{ij}\dd\bx
   =\oint_{\dOm}x_j\,\sigma_{ik}\,n_k\dd s
   =\oint_{\dOm}x_j\,(\bsig\cdot\bn)_i\dd s .
\end{equation}
By the traction-free condition \eqref{eq:tractionfree}, $\bsig\,\bn=\bm{0}$ on
all of $\dOm$, so the right-hand side of \eqref{eq:avgstress} vanishes for every
$i,j$. In particular $\int_{\Om}\sigma_{xx}\dd\bx=\int_{\Om}\sigma_{yy}\dd\bx=0$,
hence $\int_{\Om}\tr\bsig\dd\bx=0$. Substituting into \eqref{eq:dVsplit} yields
\eqref{eq:massvolume}.
\end{proof}

\noindent
Identity \eqref{eq:avgstress} is the classical statement that, in the absence of
body forces, the volume average of the stress in a body is determined entirely
by the boundary tractions; \eqref{eq:massvolume} is its specialization to a
self-equilibrated (traction-free) body, combined with the concentration
pre-stress \eqref{eq:constitutive}.

\subsection*{Comments}

\begin{remark}[The traction-free boundary is the real hypothesis]
The only geometric input to the proof is the divergence theorem; the shape of
$\Om$ plays no role. The genuine hypothesis is that the \emph{entire} boundary
is traction-free, \eqref{eq:tractionfree}. If a portion $\Gamma_{\!D}\subset\dOm$
is instead clamped ($\bu=\bm{0}$) or carries a prescribed nonzero traction
$\bm{t}$, then \eqref{eq:avgstress} retains the boundary contribution over that
portion and \eqref{eq:massvolume} acquires a correction:
\begin{equation}\label{eq:corr}
  \Delta V=\frac{\beta}{\lambda+\mu}\int_{\Om}C\dd\bx
          +\frac{1}{2(\lambda+\mu)}
            \oint_{\dOm}\bx\cdot(\bsig\cdot\bn)\dd s,
\end{equation}
where the boundary integral reduces to a sum over the parts of $\dOm$ on which
$\bsig\,\bn\neq\bm{0}$ (on a clamped part this is the reaction traction, which is
generally nonzero). Equation~\eqref{eq:massvolume} is thus exact only in the
fully stress-free setting; on the mixed boundary used in the atherosclerosis
simulations, where the lower wall $y=0$ is clamped
($u_x=u_y=0$), the exact relation is \eqref{eq:corr} and \eqref{eq:massvolume}
holds only up to the reaction-traction term produced at the clamped wall.
\end{remark}

\begin{remark}[The symmetric coupling $\beta_2=\gamma_2$ is essential]
The argument relies on the concentration entering the momentum balance and the
constitutive law through the \emph{same} coefficient, so that
\eqref{eq:equilibrium} holds for the very tensor $\bsig$ whose trace appears in
\eqref{eq:trace}. If $\beta_2\neq\gamma_2$, the momentum operator carries
$\gamma_2\grad C$ while the stress carries the pre-stress $-\beta_2 C\,\bm{I}$;
equilibrium then holds for the effective tensor
$\tilde{\bsig}=2\mu\bveps+\lambda(\divv\bu)\bm{I}-\gamma_2 C\,\bm{I}$, and the
same computation gives
\[
  \Delta V=\frac{\gamma_2}{\lambda+\mu}\int_{\Om}C\dd\bx
\]
under \eqref{eq:tractionfree} interpreted for $\tilde{\bsig}$. The single-constant
statement \eqref{eq:massvolume} is therefore special to the case
$\beta_2=\gamma_2=\beta$ adopted in the paper.
\end{remark}

\begin{remark}[Small-deformation regime]
Both ingredients are first order in the displacement: the linearized area change
$\Delta V=\int_{\Om}\divv\bu$ (the quadratic term
$\tfrac12\!\int(\partial_x u_x\,\partial_y u_y-\partial_y u_x\,\partial_x u_y)$
in $\det(\bm{I}+\grad\bu)$ is neglected) and the linear constitutive law
\eqref{eq:constitutive}. This is the same order of approximation under which the
quadratic (triangle) contribution is dropped in the rectangular derivation, so
Theorem~\ref{thm:main} holds in exactly the regime already assumed in the paper.
\end{remark}

\begin{remark}[Relation to the rectangular derivation]
For a rectangle with all four sides traction-free, \eqref{eq:avgstress}
reproduces, term by term, the boundary identities obtained in the appendix of
the paper by integrating the two elasticity equations separately; the equality
of the horizontal and vertical volume-change integrals used there is the
special case of \eqref{eq:avgstress} with $i=j$. Theorem~\ref{thm:main} may
therefore replace, or accompany, that computation: it removes the dependence on
the rectangular geometry and makes the role of the boundary condition explicit.
\end{remark}

\end{document}